\documentclass[12pt]{amsart}
\usepackage{amssymb}
\usepackage{amsmath}
\usepackage[textwidth=16cm, textheight=22cm]{geometry}
\usepackage{hyperref}

\usepackage{graphicx}
\usepackage{latexsym}
\newtheorem{algorithm}{Algorithm}
\newtheorem{example}{Example}
\newtheorem{remark}{Remark}

\def\msn{\medskip\noindent}

\def\R{\mathbb{R}}

\def\Mat{\operatorname{Mat}}

\def\bzero{O}
\def\bunity{{\mathbf 1}}

\def\semr{{\mathcal S}}

\def\msn{\smallskip\noindent}
\def\ssn{\medskip\noindent}
\def\sssn{\noindent}

\newcommand{\0}{\mathbf{0}}
\newcommand{\1}{\mathbf{1}}

\newcommand{\rseth}{\widehat{\mathbf{R}}}
\newcommand{\rset}{\mathbf{R}}
\newcommand{\rmax}{\rset_{\max}}
\newcommand{\rmin}{\rset_{\min}}
\newcommand{\rmaxh}{\rseth_{\max}}
\newcommand{\rminh}{\rseth_{\min}}
\newcommand{\smaxmin}{S_{\max,\min}}
\newcommand{\x}{\mathbf x}

\newcommand{\y}{\mathbf y}

\newcommand{\lx}{\underline{\x}}
\newcommand{\ly}{\underline{\y}}

\newcommand{\ux}{\overline{\x}}
\newcommand{\uy}{\overline{\y}}

\begin{document}

\title[Universal algorithms for Bellman equations]{Universal algorithms for solving the matrix Bellman equations over semirings}
\thanks{This work is supported by the RFBR-CRNF grant 11-01-93106, RFBR grant 12-01-00886-a and EPSRC grant RRAH15735}




\author{G.~L.~Litvinov}
\address{Grigory L. Litvinov, A.~A.~Kharkevich Institute for Information Transmission Problems,
19 B. Karetny per., GSP-4, Moscow, 127994, Russia}
\email{glitvinov@gmail.com}

\author{A.~Ya.~Rodionov}
\address{Anatoly Ya. Rodionov, Moscow Center for Continuous
  Mathematical Education,
119002, B. Vlasyevskiy per. 11, Moscow, Russia}
\email{ayarodionov@yahoo.com}

\author{S.~N.~Sergeev}
\address{Serge\u{\i} N. Sergeev, University of Birmingham,
School of Mathematics, Watson Building, Edgbaston B15 2TT, UK.}
\email{sergiej@gmail.com}

\author{A.~N.~Sobolevski}
\address{Andre\u{\i} N. Sobolevski, A.~A.~Kharkevich Institute for Information Transmission Problems,
19 B. Karetny per., GSP-4, Moscow, 127994, Russia}
\email{ansobol@gmail.com}

\keywords{idempotent semiring, tropical linear algebra, max-plus algebra, universal algorithms,
Bellman equation}

\begin{abstract}
This paper is a survey on universal algorithms for solving the matrix Bellman equations over semirings and especially tropical and idempotent semirings.
However, original algorithms are also presented. Some applications and
software implementations are discussed.
\end{abstract}

\maketitle

\section{Introduction}

Computational algorithms are constructed
on the basis of certain primitive operations. These operations manipulate
data that describe ``numbers.'' These ``numbers'' are elements of a
``numerical domain,'' that is, a mathematical object such as the field of
real numbers, the ring of integers, different semirings etc.

In practice, elements of the numerical domains are replaced
by their computer representations, that is, by elements of certain finite
models of these domains. Examples of models that can be conveniently used
for computer representation of real numbers are provided by various
modifications of floating point arithmetics, approximate arithmetics of
rational numbers~\cite{LRTch-08}, interval arithmetics etc. The difference
between mathematical objects (``ideal'' numbers) and their finite
models (computer representations) results in computational (for instance,
rounding) errors.

An algorithm is called {\it universal\/} if it is independent of a
particular numerical domain and/or its computer representation~\cite{LM-96,LM-98,LMR-00,LMa-00}.
A typical example of a universal algorithm is the computation of the
scalar product $(x,y)$ of two vectors $x=(x_1,\dots,x_n)$ and
$y=(y_1,\dots,y_n)$ by the formula $(x,y)=x_1y_1+\dots+x_ny_n$.
This algorithm (formula) is independent of a particular
domain and its computer implementation, since the formula is
well-defined for any semiring. It is clear that one algorithm can be
more universal than another. For example, the simplest Newton--Cotes formula, the
rectangular rule, provides the most universal algorithm for
numerical integration. In particular, this formula is valid also for
idempotent integration (that is, over any idempotent semiring, see
\cite{KM:97,Lit-07}).
Other quadrature formulas (for instance, combined trapezoid rule or the Simpson
formula) are independent of computer arithmetics and can be
used (for instance, in the iterative form) for computations with
arbitrary accuracy. In contrast, algorithms based on
Gauss--Jacobi formulas are designed for fixed accuracy computations:
they include constants (coefficients and nodes of these formulas)
defined with fixed accuracy. (Certainly, algorithms of this type can
be made more universal by including procedures for computing the
constants; however, this results in an unjustified complication of the
algorithms.)

Modern achievements in software development
and mathematics make us consider numerical algorithms and their
classification from a new point of view. Conventional numerical
algorithms are oriented to software (or hardware) implementation based
on floating point arithmetic and fixed accuracy. However,
it is often desirable to perform computations with variable (and
arbitrary) accuracy. For this purpose, algorithms are required
that are independent of the accuracy of computation and of the
specific computer representation of numbers. In fact, many
algorithms are independent not only of the computer representation
of numbers, but also of concrete mathematical (algebraic) operations
on data. In this case, operations themselves may be considered as variables.
Such algorithms are implemented in the form of {\it generic programs} based on
abstract data types that are
defined by the user in addition to the predefined types provided by the
language. The corresponding program tools appeared as
early as in Simula-67, but modern object-oriented languages (like
$C++$, see, for instance, \cite{Lor:93,Pohl:97}) are more convenient for generic programming. Computer algebra algorithms used in such systems as Mathematica,
Maple, REDUCE, and others are also highly universal.

A different form of universality is featured by
iterative algorithms (beginning with the successive approximation
method) for solving differential equations (for instance, methods of
Euler, Euler--Cauchy, Runge--Kutta, Adams, a number of important
versions of the difference approximation method, and the like),
methods for calculating elementary and some special functions based on
the expansion in Taylor's series and continuous fractions
(Pad\'e approximations). These algorithms are independent of the computer
representation of numbers.

The concept of a generic program was introduced by many authors;
for example, in~\cite{Leh-77} such programs were called `program schemes.'
In this paper, we discuss  universal algorithms implemented in the form of generic
programs and their specific features. This paper is closely related
to~\cite{Lit-07,LM-96,LM-98,LMa-00,LMR-00,LMRS,Ser-toep}, in which the concept of a universal algorithm was
defined and software and hardware implementation of such algorithms was
discussed in connection with problems of idempotent
mathematics, see, for instance,~\cite{KM:97,LS-01,Mik-06,Vir-00,Vir-08}.

The so-called \emph{idempotent correspondence principle}, see~\cite{LM-96,LM-98},
linking this mathematics with the usual mathematics over fields, will be
discussed below. In a nutshell,
there exists a correspondence between interesting, useful, and
important constructions and results concerning the field of real (or
complex) numbers and similar constructions dealing with various idempotent
semirings. This correspondence can be formulated in the spirit of the
well-known N.~Bohr's \emph{correspondence principle} in quantum mechanics;
in fact, the two principles are closely connected (see~\cite{Lit-07,LM-96,LM-98}).  In a sense, the traditional
mathematics over numerical fields can be treated as a `quantum' theory, whereas the idempotent mathematics
can be treated as a `classical' shadow (or counterpart) of the traditional
one. It is important that the idempotent correspondence principle is valid for algorithms, computer programs and hardware units.

In quantum mechanics the \emph{superposition principle}  means that the
Schr\"odi\-n\-ger equation (which is basic for the theory) is linear.
Similarly in idempotent mathematics the (idempotent) superposition
principle (formulated by V.~P.~Maslov)
means that some important and basic problems and equations that are nonlinear in the usual sense (for instance,
the Hamilton-Jacobi equation, which is basic for classical mechanics and appears in many
optimization problems, or the Bellman equation and its versions and
generalizations) can be treated as
linear over appropriate idempotent semirings, see~\cite{Mas-87a,Mas-87b}.

Note that numerical algorithms for infinite dimensional linear problems
over idempotent semirings (for instance, idempotent integration, integral operators
and transformations, the Hamilton--Jacobi and generalized Bellman
equations) deal with the corresponding finite-dimensional approximations.
Thus idempotent linear algebra is the basis of the idempotent numerical
analysis and, in particular, the \emph{discrete optimization theory}.

B.~A.~Carr\'e~\cite{Car-71,Car:79} (see also~\cite{Gon-75,GM:79,GM:10}) used
the idempotent linear algebra to show that different optimization problems
for finite graphs can be formulated in a unified manner and reduced to
solving Bellman equations, that is, systems of linear algebraic equations over
idempotent semirings.  He also generalized principal algorithms of
computational linear algebra to the idempotent case and showed that some of
these coincide with algorithms independently developed for solution of
optimization problems. For example, Bellman's method of
solving the shortest path problem corresponds to a version of Jacobi's
method for solving a system of linear equations, whereas Ford's algorithm
corresponds to a version of Gauss--Seidel's method. We briefly discuss Bellman equations and the corresponding optimization problems on graphs, and use the
ideas of Carr\'{e} to obtain new universal algorithms. We stress that these well-known results can be interpreted as a
manifestation of the idempotent superposition principle.

Note that many algorithms for solving the matrix Bellman equation could be found in~\cite{BCOQ,Car-71,Car:79,CG:79,Gon-75,GM:10,LMRS,LMa-00,Rot-85,Ser-toep,CS-07}.
More general problems of linear algebra over the max-plus algebra are examined, for instance in~\cite{But:10}.

We also briefly discuss interval analysis over idempotent and positive semirings. Idempotent interval analysis appears in~\cite{LS-00,LS-01,Sob-99}, where it
is applied to the Bellman matrix equation. Many different problems coming from the idempotent linear algebra, have been considered since then, see for instance~\cite{CC-02,Fie+06,Har+09,Mys-05,Mys-06}.
It is important to observe that intervals over an idempotent semiring form a new idempotent semiring. Hence universal algorithms can be applied to elements of this new semiring and generate interval extensions of the initial algorithms.

This paper is about software implementations of universal algorithms for solving the matrix Bellman equations over semirings. In Section~\ref{s:sem} we pre\-sent an introduction to mathematics of semirings and especially to
the tropical (idempotent) mathematics, that is, the area of mathematics
working with  \emph{idempotent semirings} (that is, semirings with idempotent addition). In Section~\ref{s:main} we present a number of
well-known and new universal algorithms of linear algebra over semirings, related to discrete matrix Bellman equation and algebraic path problem.
These algorithms are closely related to their linear-algebraic prototypes
described, for instance, in the celebrated book of Golub and Van Loan~\cite{GvL} which serves as the main source of such prototypes. Following the
style of~\cite{GvL} we pre\-sent them in MATLAB code. The perspectives and experience of their implementation are also discussed.

\if{

\section{Universal algorithms and accuracy of computations}

Calculations on computers usually are based on a floating-point arithmetic
with a mantissa of a fixed length; that is, computations are performed
with fixed accuracy. Broadly speaking, with this approach only
the relative rounding error is fixed, which can lead to a drastic
loss of accuracy and invalid results (for instance, when summing series and
subtracting close numbers). On the other hand, this approach provides
rather high speed of computations. Many important numerical algorithms
are designed to use floating-point arithmetic (with fixed accuracy)
and ensure the maximum computation speed. However, these algorithms
are not universal. The above mentioned Gauss--Jacobi quadrature formulas,
computation of elementary and special functions on the basis of the
best polynomial or rational approximations or Pad\'e--Chebyshev
approximations, and some others belong to this type. Such algorithms
use nontrivial constants specified with fixed accuracy.

Recently, problems of accuracy, reliability, and authenticity of
computations (including the effect of rounding errors) have gained
much atention; in part, this fact is related to the ever-increasing
performance of computer hardware. When errors in initial data and
rounding errors strongly affect the computation results, such as in ill-posed
problems, analysis of stability of solutions, etc., it is often useful
to perform computations with improved and variable accuracy. In
particular, the rational arithmetic, in which the rounding error is
specified by the user [44], can be used for this purpose.
This arithmetic is a useful complement to the interval analysis [54].
The corresponding computational algorithms must be
universal (in the sense that they must be independent of the computer
representation of numbers).

}\fi

\section{Mathematics of semirings}
\label{s:sem}

\subsection{Basic definitions}

A broad class of universal algorithms is related to the concept of
a semiring. We recall here the definition (see, for instance, \cite{Gol:99}).

A set~$S$ is called a \emph{semiring} if it is endowed with two
associative operations: \emph{addition}~$\oplus$ and
\emph{multiplication}~$\odot$ such that addition is commutative,
multiplication distributes over addition from either side, $\0$ (resp.,
$\1$) is the neutral element of addition (resp., multiplication), $\0 \odot
x = x \odot \0 = \0$ for all $x \in S$, and $\0 \neq \1$.

Let the semiring~$S$ be partially ordered by a relation~$\preceq$ such that $\0$ is
the least element and the inequality $x \preceq y$ implies that $x \oplus z
\preceq y \oplus z$, $x \odot z \preceq y \odot z$, and~$z \odot x \preceq z \odot
y$ for all $x, y, z \in S$; in this case the semiring~$S$ is called
\emph{positive} (see, for instance, \cite{Gol:99}).

An element $x\in S$ is called \emph{invertible} if there exists an
element $x^{-1}\in S$ such that $xx^{-1}=x^{-1}x=\1$.
A semiring~$S$ is called a \emph{semifield} if every nonzero element is invertible.

A semiring~$S$ is called \emph{idempotent} if $x \oplus x = x$ for all
$x \in S$. 
In this case the
addition~$\oplus$ defines a \emph{canonical
partial order}~$\preceq$ on the semiring~$S$ by the rule: $x\preceq y$ iff $x\oplus y=y$. It is easy to prove that any idempotent semiring is
positive with respect to this order. Note also that $x \oplus y =
\sup\{x,y\}$ with respect to the canonical order. In the sequel, we shall
assume that all idempotent semirings are ordered by the canonical partial
order relation.

We shall say that a positive (for instance, idempotent) semiring $S$ is \emph{complete} if  for every subset $T\subset S$ there exist elements $\sup T\in S$ and $\inf T\in S$, and if the operations $\oplus$ and $\odot$ distribute over such sups and infs.

The most well-known and important examples of positive semirings are
``numerical'' semirings consisting of (a subset of) real numbers and ordered
by the usual linear order $\leq$ on~$\rset$: the semiring~$\rset_+$
with the usual operations $\oplus = +$, $\odot = \cdot$ and neutral
elements $\0 = 0$, $\1 = 1$, the semiring~$\rmax = \rset \cup \{-\infty\}$
with the operations $\oplus = \max$, $\odot = +$ and neutral elements $\0 =
-\infty$, $\1 = 0$, the semiring $\rmaxh = \rmax \cup \{\infty\}$,
where $x \preceq \infty$, $x \oplus \infty = \infty$ for all $x$, $x \odot
\infty = \infty \odot x = \infty$ if $x \neq \0$, and $\0 \odot \infty =
\infty \odot \0$, and the semiring~$\smaxmin^{[a,b]} = [a, b]$, where
$-\infty \leq a < b \leq +\infty$, with the operations $\oplus =
\max$, $\odot = \min$ and neutral elements $\0 = a$, $\1 = b$.  The
semirings~$\rmax$, $\rmaxh$, and~$\smaxmin^{[a,b]} = [a, b]$ are
idempotent. The semirings $\rmaxh$, $S^{[a,b]}_{\mathrm{max, min}}$, $\widehat{\textbf{R}}_+=\textbf{R}_+\bigcup\{\infty\}$ are complete.
Remind that every partially ordered set can be imbedded to its completion (a minimal complete set containing the initial one).
The semiring $\textbf{R}_\mathrm{min}=\textbf{R}\bigcup\{\infty\}$ with operations $\oplus=\mathrm{min}$ and $\odot=+$ and neutral elements $\textbf{0}=\infty$, $\textbf{1}=0$ is isomorphic to $\rmax$.

The semiring $\rmax$ is also called the \emph{max-plus algebra}. The semifields $\rmax$ and $\rmin$ are called \emph{tropical algebras}.
The term ``tropical'' initially appeared in~\cite{Sim-88}
for a discrete version of the max-plus algebra as a suggestion of Ch.~Choffrut, see also~\cite{Gun:98,Mik-06,Vir-08}.

Many mathematical constructions, notions, and results over the fields of
real and complex numbers have nontrivial analogs over idempotent semirings.
Idempotent semirings have become recently the object of investigation of
new branches of mathematics, \emph{idempotent mathematics} and \emph{tropical geometry}, see, for instance~\cite{BCOQ,CG:79,Lit-07,Mik-06,Vir-00,Vir-08}.

Denote by $\Mat_{mn}(S)$ a set of all matrices $A = (a_{ij})$
with $m$~rows and $n$~columns whose coefficients belong to a semiring~$S$.
The sum $A \oplus B$ of matrices $A, B \in \Mat_{mn}(S)$ and the product
$AB$ of matrices $A \in \Mat_{lm}(S)$ and $B \in \Mat_{mn}(S)$ are defined
according to the usual rules of linear algebra:
$A\oplus B=(a_{ij} \oplus b_{ij})\in \mathrm{Mat}_{mn}(S)$ and
$$
AB=\left(\bigoplus_{k=1}^m a_{ij}\odot b_{kj}\right)\in\Mat_{ln}(S),
$$
where $A\in \Mat_{lm}(S)$ and $B\in\Mat_{mn}(S)$.
Note that we write $AB$ instead of $A\odot B$.

If the semiring~$S$ is
positive, then the set \\ $\Mat_{mn}(S)$ is ordered by the relation $A =
(a_{ij}) \preceq B = (b_{ij})$ iff $a_{ij} \preceq b_{ij}$ in~$S$ for all $1
\leq i \leq m$, $1 \leq j \leq n$.

The matrix multiplication is consistent with the order~$\preceq$ in the
following sense: if $A, A' \in \Mat_{lm}(S)$, $B, B' \in \Mat_{mn}(S)$ and
$A \preceq A'$, $B \preceq B'$, then $AB \preceq A'B'$ in $\Mat_{ln}(S)$. The set
$\Mat_{nn}(S)$ of square $(n \times n)$ matrices over a [positive,
idempotent] semiring~$S$ forms a [positive, idempotent] semi-ring with a
zero element $O = (o_{ij})$, where $o_{ij} = \0$, $1 \leq i, j
\leq n$, and a unit element $I = (\delta_{ij})$, where $\delta_{ij} =
\1$ if $i = j$ and $\delta_{ij} = \0$ otherwise.

The set $\Mat_{nn}$ is an example of a noncommutative semiring if $n>1$.


\subsection{Closure operation}

In what follows, we are mostly interested in complete positive semirings, and particularly
in idempotent semirings. Regarding examples of the previous section, recall that the semirings  $S^{[a,b]}_{\mathrm{max, min}}$, $\rmaxh=\rmax\cup\{+\infty\}$,
$\rminh=\rmin\cup\{-\infty\}$
and $\widehat{\textbf{R}}_+=\textbf{R}_+\cup\{+\infty\}$ are complete positive, and the semirings
$S^{[a,b]}_{\mathrm{max, min}}$, $\rmaxh$ and $\rminh$ are idempotent.

$\widehat{\textbf{R}}_+$ is a completion of $\textbf{R}_+$, and
$\rmaxh$ (resp. $\rminh$) are completions of $\rmax$ (resp.
$\rmin$). More generally, we note that any positive semifield $S$
can be completed by means of a standard procedure, which uses
Dedekind cuts and is described in~\cite{Gol:99,LMSz-01}. The result
of this completion is a semiring $\widehat{S}$, which is not a
semifield anymore.

The semiring of matrices $\Mat_{nn}(S)$ over a complete positive (resp., idempotent) semiring is again a complete positive (resp., idempotent) semiring. For more background in complete idempotent semirings, the
reader is referred to~\cite{LMSz-01}.

In any complete positive semiring $S$ we have a unary operation of
{\em closure} $a\mapsto a^*$ defined by
\begin{equation}
\label{cldef-gen}
a^*:=\sup_{k\geq 0} \1\oplus a\oplus \ldots\oplus a^k,
\end{equation}
Using that the operations $\oplus$ and $\odot$ distribute over such sups,
it can be shown that $a^*$ is the {\bf least solution} of $x=ax\oplus\1$ and $x=xa\oplus\1$, and also that $a^*b$ is the
the least solution
of $x=ax\oplus b$ and $x=xa\oplus b$.

In the case of idempotent addition~\eqref{cldef-gen}
becomes particularly nice:
\begin{equation}
\label{aldef-id}
a^*=\bigoplus_{i\geq 0} a^i=\sup_{i\geq 0} a^i.
\end{equation}

If a positive semiring $S$ is not complete, then it often happens that the closure operation
can still be defined on some ``essential'' subset of $S$. Also recall that
any positive semifield $S$ can be completed~\cite{Gol:99,LMSz-01}, and then
the closure is defined for every element of the completion.

In numerical semirings the operation~$*$ is usually very easy to implement:
$x^* = (1-x)^{-1}$ if $x<1$ in $\rset_+$, or $\widehat{\textbf{R}}_+$ and $x^*=\infty$ if $x\geq 1$ in $\widehat{\textbf{R}}_+$; $x^* = \1$ if $x \preceq \1$ in $\rmax$
and $\rmaxh$, $x^* = \infty$ if $x \succ \1$ in $\rmaxh$, $x^* = \1$
for all $x$ in $\smaxmin^{[a,b]}$. In all other cases $x^*$ is undefined.

The closure operation in matrix semirings over a complete positive semiring~$S$
can be defined as in~\eqref{cldef-gen}:
\begin{equation}
\label{cldef-matx}
A^*:=\sup_{k\geq 0} I\oplus A\oplus \ldots\oplus A^k,
\end{equation}
and one can show that it is the least solution $X$ satisfying the matrix equations
$X=AX\oplus I$ and $X=XA\oplus I$.

Equivalently, $A^*$ can
be defined by induction: let $A^*
= (a_{11})^* = (a^*_{11})$ in $\Mat_{11}(S)$ be defined by~\eqref{cldef-gen},
and for any integer $n > 1$
and any matrix
$$
   A = \begin{pmatrix} A_{11}& A_{12}\\ A_{21}& A_{22} \end{pmatrix},
$$
where $A_{11} \in \Mat_{kk}(S)$, $A_{12} \in \Mat_{k\, n - k}(S)$,
$A_{21} \in \Mat_{n - k\, k}(S)$, $A_{22} \in \Mat_{n - k\, n - k}(S)$,
$1 \leq k \leq n$, by definition,
\begin{equation}
\label{A_Star}
   A^* = \begin{pmatrix}
   A^*_{11} \oplus A^*_{11} A_{12} D^* A_{21} A^*_{11} &
   \quad A^*_{11} A_{12} D^* \\[2ex]
   D^* A_{21} A^*_{11} &
   D^*
   \end{pmatrix},
\end{equation}
where $D = A_{22} \oplus A_{21} A^*_{11} A_{12}$.

Defined here for complete positive semirings, the closure operation is a semiring
analogue of the operation $(1-a)^{-1}$ and, further, $(I-A)^{-1}$ in matrix algebra over the
field of real or complex mumbers.
This operation can be thought of as {\bf regularized sum} of the
series $I+A+A^2+\ldots$, and the closure operation defined above is another such
regularization. Thus we can also define the closure operation $a^*=(1-a)^{-1}$
and $A^*=(I-A)^{-1}$ in the traditional linear algebra.
To this end, note that the recurrence
relation above coincides with the formulas of escalator method of matrix
inversion in the traditional linear algebra over the field of real or
complex numbers, up to the algebraic operations used. Hence this algorithm
of matrix closure requires a polynomial number of operations in~$n$, see
below for more details.

Let~$S$ be a complete positive semiring. The \emph{matrix (or discrete stationary)
Bellman equation} has the form
\begin{equation}
\label{AX}
    X = AX \oplus B,
\end{equation}
where $A \in \Mat_{nn}(S)$, $X, B \in \Mat_{ns}(S)$, and the matrix~$X$ is
unknown. As in the scalar case, it
can be shown that for complete positive semirings, if $A^*$ is defined as in~\eqref{cldef-matx} then $A^*B$ is the least in the set of solutions to equation~\eqref{AX}
with respect to the partial order in $\Mat_{ns}(S)$. In the
idempotent case
\begin{equation}
\label{aldef-matx}
A^*=\bigoplus_{i\geq 0} A^i=\sup_{i\geq 0} A^i.
\end{equation}

Consider also the case when $A=(a_{ij})$ is $n\times n$ {\em strictly upper-triangular}
(such that $a_{ij}=\0$ for $i\geq j$), or $n\times n$ {\em strictly lower-triangular}
(such that $a_{ij}=\0$ for $i\leq j$). In this case $A^n=O$, the all-zeros matrix,
and it can be shown by iterating $X=AX\oplus I$ that this equation has a unique
solution, namely
\begin{equation}
\label{a*nilp}
A^*=I\oplus A\oplus\ldots\oplus A^{n-1}.
\end{equation}
Curiously enough, formula~\eqref{a*nilp} works more generally in the
case of numerical idempotent semirings:
in fact, the series~\eqref{aldef-matx}
converges there if and only if it can be truncated to~\eqref{a*nilp}. This
is closely related to the principal path interpretation of $A^*$ explained in
the next subsection.

In fact, theory of the discrete stationary Bellman equation can be
developed using the identity $A^* = AA^* \oplus I$ as an axiom
without any explicit formula for the closure (the so-called \emph{closed
semirings}, see, for instance, \cite{Gol:99,Leh-77,Rot-85}). Such theory can be based
on the following identities, true both for the case of idempotent
semirings and the real numbers with
conventional arithmetic (assumed that $A$ and $B$ have appropriate sizes):
\begin{equation}
\label{e:conway}
\begin{split}
(A\oplus B)^*&=(A^*B)^*A^*,\\
(AB)^*A&=A(BA)^*.
\end{split}
\end{equation}
This abstract setting unites the case of positive and idempotent semirings with
the conventional linear algebra over the field of real and complex numbers.

\subsection{Weighted directed graphs and matrices over semirings}

Suppose that $S$ is a semiring with zero~$\0$ and unity~$\1$. It is well-known
that any square matrix $A = (a_{ij}) \in \Mat_{nn}(S)$ specifies a
\emph{weighted directed graph}. This geometrical construction includes
three kinds of objects: the set $X$ of $n$ elements $x_1, \dots, x_n$
called \emph{nodes}, the set $\Gamma$ of all ordered pairs $(x_i, x_j)$
such that $a_{ij} \neq \0$ called \emph{arcs}, and the mapping $A \colon
\Gamma \to S$ such that $A(x_i, x_j) = a_{ij}$. The elements $a_{ij}$ of
the semiring $S$ are called \emph{weights} of the arcs.Conversely, any given weighted
directed graph with $n$ nodes specifies a
unique matrix $A \in \Mat_{nn}(S)$.

\begin{figure}[h]
\centering
\setlength{\belowcaptionskip}{25pt}
\caption[]{}
\label{Fig1.}
\includegraphics[scale=0.55]{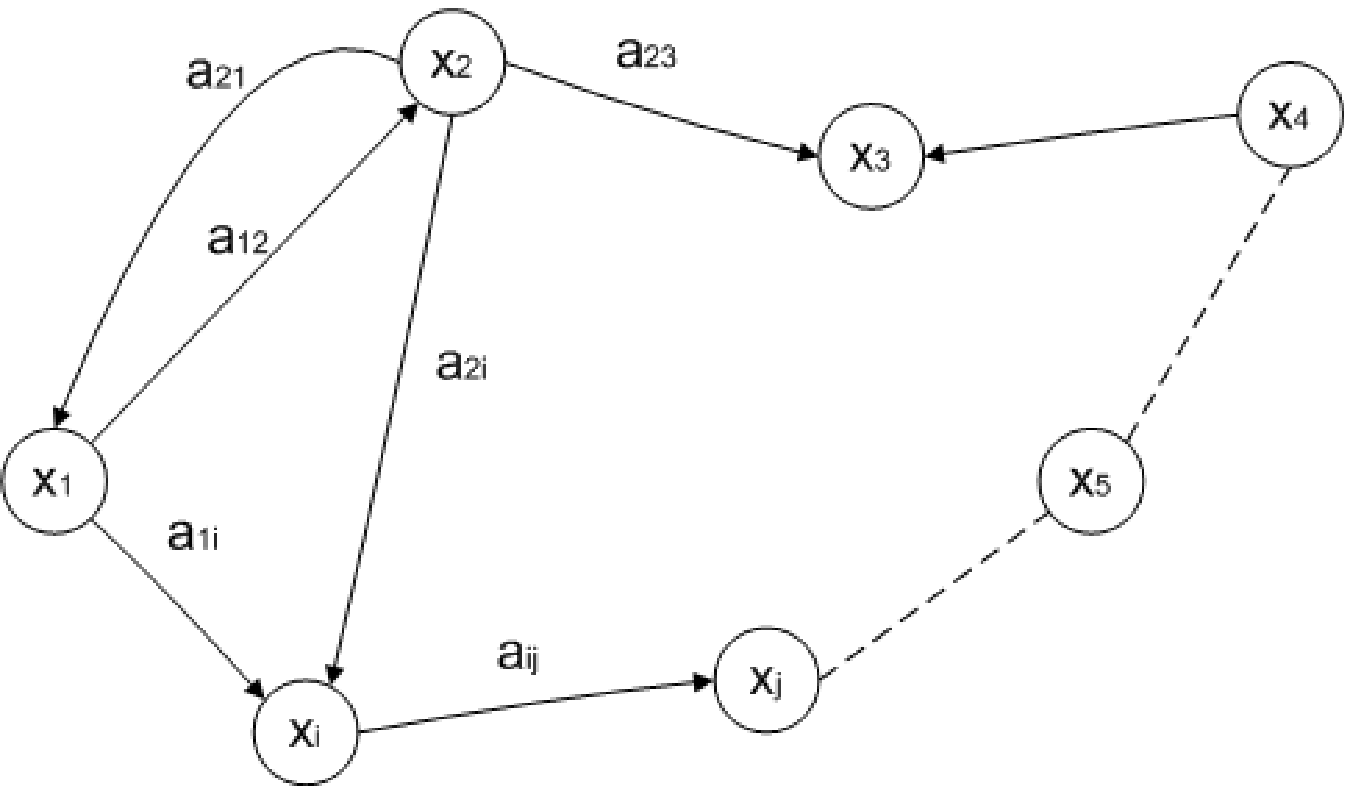}
\end{figure}

This definition allows for some pairs of nodes to be disconnected if the
corresponding element of the matrix $A$ is $\0$ and for some channels to be
``loops'' with coincident ends if the matrix $A$ has nonzero diagonal
elements.

Recall that a sequence of nodes of the form
$$
    p = (y_0, y_1, \dots, y_k)
$$
with $k \geq 0$ and $(y_i, y_{i + 1}) \in \Gamma$, $i = 0, \dots, k -
1$, is called a \emph{path} of length $k$ connecting $y_0$ with $y_k$.
Denote the set of all such paths by $P_k(y_0,y_k)$. The weight $A(p)$ of a
path $p \in P_k(y_0,y_k)$ is defined to be the product of weights of arcs
connecting consecutive nodes of the path:
$$
    A(p) = A(y_0,y_1) \odot \cdots \odot A(y_{k - 1},y_k).
$$
By definition, for a `path' $p \in P_0(x_i,x_j)$ of length $k = 0$ the
weight is $\1$ if $i = j$ and $\0$ otherwise.

For each matrix $A \in \Mat_{nn}(S)$ define $A^0 = E = (\delta_{ij})$
(where $\delta_{ij} = \1$ if $i = j$ and $\delta_{ij} = \0$ otherwise) and
$A^k = AA^{k - 1}$, $k \geq 1$.  Let $a^{[k]}_{ij}$ be the $(i,j)$th
element of the matrix $A^k$. It is easily checked that
$$
   a^{[k]}_{ij} =
   \bigoplus_{\substack{i_0 = i,\, i_k = j\\
    1 \leq i_1, \ldots, i_{k - 1} \leq n}}
    a_{i_0i_1} \odot \dots \odot a_{i_{k - 1}i_k}.
$$
Thus $a^{[k]}_{ij}$ is the supremum of the set of weights corresponding to
all paths of length $k$ connecting the node $x_{i_0} = x_i$ with $x_{i_k} =
x_j$.

Let $A^*$ be defined as in~\eqref{aldef-matx}.
Denote the elements of the matrix $A^*$ by $a^{*}_{ij}$, $i, j = 1,
\dots, n$; then
$$
    a^{*}_{ij}
    = \bigoplus_{0 \leq k < \infty}
    \bigoplus_{p \in P_k(x_i, x_j)} A(p).
$$

The closure matrix $A^*$ solves the well-known \emph{algebraic path
problem}, which is formulated as follows: for each pair $(x_i,x_j)$
calculate the supremum of weights of all paths (of arbitrary length)
connecting node $x_i$ with node $x_j$. The closure operation in matrix
semirings has been studied extensively (see, for instance,
\cite{BCOQ,Car-71,Car:79,CG:79,Gol:99,GM:79,GM:10,KM:97,LS-01} and references therein).

\begin{example}[The shortest path problem]
{\rm Let $S = \rmin$, so the weights are real numbers. In this case
$$
    A(p) = A(y_0,y_1) + A(y_1,y_2) + \dots + A(y_{k - 1},y_k).
$$
If the element $a_{ij}$ specifies the length of the arc $(x_i,x_j)$ in some
metric, then $a^{*}_{ij}$ is the length of the shortest path connecting
$x_i$ with $x_j$.}
\end{example}

\begin{example}[The maximal path width problem]
{\rm Let $S = \rset \cup \{\0,\1\}$ with $\oplus = \max$, $\odot = \min$. Then
$$
    a^{*}_{ij} =
    \max_{p \in \bigcup\limits_{k \geq 1} P_k(x_i,x_j)} A(p),
$$
$$
        A(p) = \min (A(y_0,y_1), \dots, A(y_{k - 1},y_k)).
$$
If the element $a_{ij}$ specifies the ``width'' of the arc
$(x_i,x_j)$, then the width of a path $p$ is defined as the minimal
width of its constituting arcs and the element $a^{*}_{ij}$ gives the
supremum of possible widths of all paths connecting $x_i$ with $x_j$.}
\end{example}

\begin{example}[A simple dynamic programming problem]
{\rm Let $S = \rmax$ and suppose $a_{ij}$ gives the \emph{profit} corresponding
to the transition from $x_i$ to $x_j$. Define the vector $B  = (b_i) \in
\Mat_{n1}(\rmax)$ whose element $b_i$ gives the \emph{terminal profit}
corresponding to exiting from the graph through the node $x_i$. Of course,
negative profits (or, rather, losses) are allowed. Let $m$ be the total
profit corresponding to a path $p \in P_k(x_i,x_j)$, that is
$$
    m = A(p) + b_j.
$$
Then it is easy to check that the supremum of profits that can be achieved
on paths of length $k$ beginning at the node $x_i$ is equal to $(A^kB)_i$
and the supremum of profits achievable without a restriction on the length
of a path equals $(A^*B)_i$.}
\end{example}

\begin{example}[The matrix inversion problem]
{\rm Note that in the formulas of this section we are using distributivity of
the multiplication $\odot$ with respect to the addition $\oplus$ but do not
use the idempotency axiom. Thus the algebraic path problem can be posed for
a nonidempotent semiring $S$ as well (see, for instance, \cite{Rot-85}). For
instance, if $S = \rset$, then
$$
    A^* = I + A + A^2 + \dotsb = (I - A)^{-1}.
$$
If $\|A\| > 1$ but the matrix $I - A$ is invertible, then this expression
defines a regularized sum of the divergent matrix power series
$\sum_{i \geq 0} A^i$.}
\end{example}

We emphasize that this connection between the matrix closure operation and
solutions to the Bellman equation gives rise to a number of different
algorithms for numerical calculation of the matrix closure. All these
algorithms are adaptations of the well-known algorithms of the traditional
computational linear algebra, such as the Gauss--Jordan elimination, various
iterative and escalator schemes, etc. This is a special case of the idempotent superposition principle (see below).

\subsection{Interval analysis over positive semirings}

Traditional interval analysis is a nontrivial and popular mathematical area, see, for instance,~\cite{AH:83,Fie+06,Kre+98,Moo:79,Neu:90}. An ``idempotent'' version of interval analysis (and moreover interval analysis over positive semirings) appeared in~\cite{LS-00,LS-01,Sob-99}. Rather many publications on the subject appeared later, see, for instance,~\cite{CC-02,Fie+06,Har+09,Mys-05,Mys-06}. Interval analysis over the positive semiring $\textbf{R}_+$ was discussed in~\cite{BN-74}.

Let a set~$S$ be partially ordered by a relation $\preceq$.
A \emph{closed interval} in~$S$ is a subset of the form $\x = [\lx, \ux] =
\{\, x \in S \mid \lx \preceq x \preceq \ux\, \}$, where the elements $\lx \preceq
\ux$ are called \emph{lower} and \emph{upper bounds} of the interval $\x$.
The order~$\preceq$ induces a partial ordering on the set of all closed
intervals in~$S$: $\x \preceq \y$ iff $\lx \preceq \ly$ and $\ux \preceq \uy$.

A \emph{weak interval extension} $I(S)$ of a positive semiring~$S$ is the
set of all closed intervals in~$S$ endowed with operations $\oplus$
and~$\odot$ defined as ${\x \oplus \y} = [{\lx \oplus \ly}, {\ux \oplus
\uy}]$, ${\x \odot \y} = [{\lx \odot \ly}, {\ux \odot \uy}]$ and a partial
order induced by the order in $S$. The closure operation in $I(S)$ is
defined by $\x^* = [\lx^*, \ux^*]$. There are some other interval extensions (including the so-called strong interval extension~\cite{LS-01}) but the weak extension is more convenient.

The extension $I(S)$ is positive; $I(S)$ is idempotent if $S$ is an idempotent semiring.
A universal algorithm over $S$ can be applied to $I(S)$ and we shall get an interval version of the initial algorithm.
Usually both versions have the same complexity. For the discrete stationary Bellman equation and the corresponding optimization problems on graphs, interval analysis was examined in~\cite{LS-00,LS-01} in details. Other problems of idempotent linear algebra were examined in~\cite{CC-02,Fie+06,Har+09,Mys-05,Mys-06}.

Idempotent mathematics appears to be remarkably simpler than its
traditional analog. For example, in traditional interval arithmetic,
multiplication of intervals is not distributive with respect to addition of
intervals, whereas in idempotent interval arithmetic this distributivity is
preserved. Moreover, in traditional interval analysis the set of all
square interval matrices of a given order does not form even a semigroup
with respect to matrix multiplication: this operation is not associative
since distributivity is lost in the traditional interval arithmetic. On the
contrary, in the idempotent (and positive) case associativity is preserved. Finally, in
traditional interval analysis some problems of linear algebra, such as
solution of a linear system of interval equations, can be very difficult
(more precisely, they are $NP$-hard, see~
\cite{Kre+98} and references therein). It was noticed  in~\cite{LS-00,LS-01} that in the idempotent case solving an interval linear system
requires a polynomial number of operations (similarly to the usual Gauss
elimination algorithm).  Two properties that make the idempotent interval
arithmetic so simple are monotonicity of arithmetic operations and
positivity of all elements of an idempotent semiring.

Interval estimates in idempotent mathematics are usually exact.
In the traditional theory such estimates tend to be overly pessimistic.

\subsection{Idempotent correspondence principle}

There is a nontrivial analogy between mathematics of semirings and
quantum mechanics. For example, the field of real numbers can be
treated as a ``quantum object'' with respect to idempotent semirings.
So idempotent semirings can be treated as ``classical'' or
``semi-classical'' objects with respect to the field of real
numbers.

Let $\rset$ be the field of real numbers and $\rset_+$ the subset
of all non-negative numbers. Consider the following change of variables:
$$
u \mapsto w = h \ln u,
$$
where $u \in \rset_+ \setminus \{0\}$, $h > 0$; thus $u = e^{w/h}$,
$w \in
\rset$. Denote by $\0$ the additional element $-\infty$ and by $S$
the extended real line $\rset \cup \{\0\}$. The above change of
variables has a
natural extension $D_h$ to the whole $S$ by $D_h(0) = \0$; also, we
denote $D_h(1) = 0 = \1$.

Denote by $S_h$ the set $S$ equipped with the two operations $\oplus_h$
(generalized addition) and $\odot_h$ (generalized multiplication)
such that
$D_h$ is a homomorphism of $\{\rset_+, +, \cdot\}$ to $\{S, \oplus_h,
\odot_h\}$. This means that $D_h(u_1 + u_2) = D_h(u_1) \oplus_h D_h(u_2)$
and $D_h(u_1 \cdot u_2) = D_h(u_1) \odot_h D_h(u_2)$, that is, $w_1 \odot_h
w_2 = w_1 + w_2$ and $w_1 \oplus_h w_2 = h \ln (e^{w_1/h} + e^{w_2/h})$.
It is easy to prove that $w_1 \oplus_h w_2 \to \max\{w_1, w_2\}$
as $h \to 0$.

\if{
Denote by $\rmax$ the set $S = \rset \cup \{\0\}$ equipped
with operations
$\oplus = \max$ and $\odot = +$, where ${\0} = -\infty$, ${\1} = 0$
as above.
}\fi
$\rset_+$ and $S_h$ are isomorphic semirings; therefore
we have obtained $\rmax$ as a result of a deformation of $\rset_+$.
We stress the obvious analogy with the quantization procedure,
where $h$ is
the analog of the Planck constant. In these terms, $\rset_+$
(or $\rset$)
plays the part of a ``quantum object'' while $\rmax$ acts as a
``classical'' or ``semi-classical'' object that arises as the result
of a {\it dequantization} of this quantum object.
In the case of $\rmin$, the corresponding dequantization procedure is
generated by the change of variables $u \mapsto w = -h \ln u$.

There is a natural transition from the field of real numbers or complex numbers to the idempotent semiring $\rmax$ (or $\rmin$). This is a composition of the mapping $x\mapsto |x|$ and the deformation described above.

In general an \emph{idempotent dequantization} is a transition from a basic field to an idempotent semiring in mathematical concepts, constructions and results, see~\cite{Lit-07,LM-98} for details. Idempotent dequantization suggests the following formulation of the idempotent
correspondence principle:
\begin{quote}
{\it There exists a heuristic correspondence between interesting, useful,
and important constructions and results over the field of real (or
complex) numbers and similar constructions and results over idempotent
semirings in the spirit of N.~Bohr's correspondence principle in
quantum mechanics.}
\end{quote}

\begin{figure}[h]
\centering
\setlength{\belowcaptionskip}{25pt}
\label{Fig.2}
\includegraphics[scale=0.67]{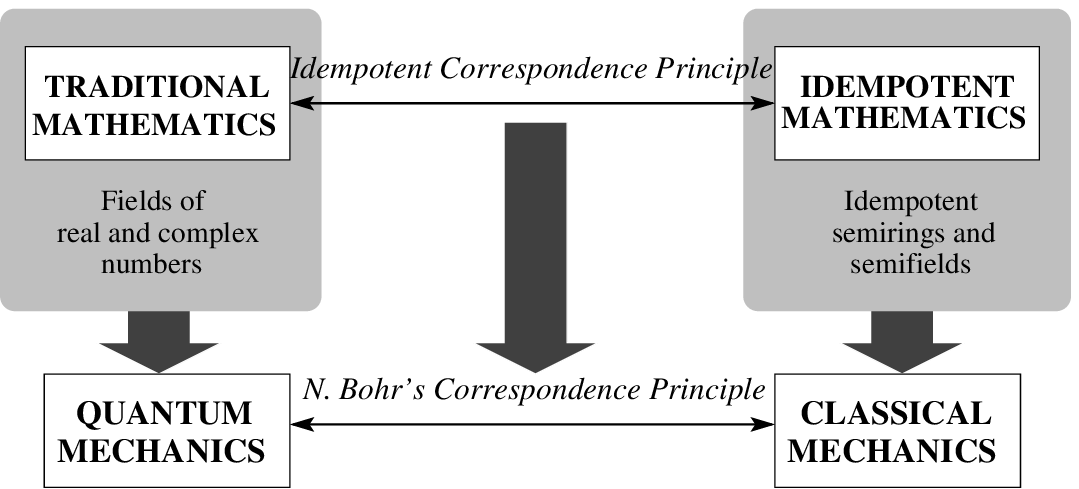}
\end{figure}

Thus idempotent mathematics can be treated as a ``classical shadow (or
counterpart)'' of the traditional Mathematics over fields.
A systematic application of this correspondence principle leads to a variety of theoretical and applied results, see,
for instance, ~\cite{Lit-07,LM-98,LMa-00,LS-01,Mik-06,Vir-00,Vir-08}.
Relations to quantum physics are discussed in detail, for instance, in~\cite{Lit-07}.

In this paper we aim to develop a practical systematic application of the correspondence principle to the algorithms of linear algebra and discrete
mathematics. For the remainder of this subsection
let us focus on an idea how the idempotent correspondence principle
may lead to a unifying approach to hardware
design. (See~\cite{LMR-00,LMRS} for more information.)

The most important and standard numerical algorithms have many hardware
realizations in the form of technical devices or special processors.
These devices often can be used as prototypes for new hardware
units resulting from mere substitution of the usual arithmetic operations
by their semiring analogs (and additional tools for generating neutral
elements $\0$ and $\1$). Of course,
the case of numerical semirings consisting of real numbers (maybe except
neutral elements)  and semirings of numerical intervals is the most simple and natural. Note that for semifields (including $\rmax$ and $\rmin$)
the operation of division is also defined.

Good and efficient technical ideas and decisions can be taken
from prototypes to new hardware units. Thus the correspondence
principle generates a regular heuristic method for hardware design.
Note that to get a patent it is necessary to present the so-called
`invention formula', that is to indicate a prototype for the suggested
device and the difference between these devices.

Consider (as a typical example) the most popular and important algorithm
of computing the scalar product of two vectors:
\begin{equation}
\label{usscalprod}
(x,y)=x_1y_1+x_2y_2+\cdots + x_ny_n.
\end{equation}
The universal version of~\eqref{usscalprod} for any semiring $A$ is obvious:
\begin{equation}
\label{semscalprod}
(x,y)=(x_1\odot y_1)\oplus(x_2\odot y_2)\oplus\cdots\oplus
(x_n\odot y_n).
\end{equation}
In the case $A=\rmax$ this formula turns into the following one:
\begin{equation}
\label{maxscalprod}
(x,y)=\max\{ x_1+y_1,x_2+y_2, \cdots, x_n+y_n\}.
\end{equation}

This calculation is standard for many optimization algorithms, so
it is useful to construct a hardware unit for computing~\eqref{maxscalprod}. There
are many different devices (and patents) for computing~\eqref{usscalprod} and every
such device can be used as a prototype to construct a new device for
computing~\eqref{maxscalprod} and even~\eqref{semscalprod}. Many processors for matrix multiplication
and for other algorithms of linear algebra are based on computing
scalar products and on the corresponding ``elementary'' devices.
Using modern technologies
it is possible to construct cheap special-purpose multi-processor
chips and systolic arrays of elementary processors implementing
universal algorithms. See, for instance, ~\cite{LMR-00,LMRS,Rot-85} where
the systolic arrays and parallel computing issues are discussed for
the algebraic path problem. In particular, there is a systolic array of
$n(n+1)$ elementary processors which performs computations of the Gauss--Jordan
elimination algorithm and can solve the algebraic path problem within $5n-2$
time steps.

\section{Some universal algorithms of linear algebra}
\label{s:main}

In this section we discuss universal algorithms computing $A^*$
and $A^*B$. We start with the basic escalator and Gauss-Jordan elimination  techniques in Subsect.~\ref{ss:GJ} and continue with its specification
to the case of Toeplitz systems in Subsect.~\ref{ss:Toep}. The universal
LDM decomposition of Bellman equations is explained in Subsect.~\ref{ss:LDM},
followed by its adaptations to symmetric and band matrices in Subsect.~\ref{ss:LDMspec}. The iteration schemes are discussed in Subsect.~\ref{ss:Iter}. In the final Subsect.~\ref{ss:impl} we discuss the implementations
of universal algorithms.

Algorithms themselves will be described in a language of Matlab,
following the tradition of Golub and van Loan~\cite{GvL}. This is done for
two purposes: 1) to simplify the comparison of the algorithms with their
prototypes taken mostly from~\cite{GvL}, 2) since the language of Matlab
is designed for matrix computations. We will not formally describe the rules
of our Matlab-derived language, preferring just to outline the following important features:
\begin{itemize}
\item[1.] Our basic arithmetic operations are $a\oplus b$, $a\odot b$ and $a^*$.
\item[2.] The vectorization of these operations follows the rules of Matlab.
\item[3.] We use basic keywords of Matlab like `for', `while', 'if' and 'end',
similar to other programming languages like $C++$ or Java.
\end{itemize}
Let us give some examples of universal matrix computations in our language:\\
{\em Example 1.}  $v(1: j)=\alpha^*\odot a(1: j,k)$ means that the result of (scalar) multiplication of the first $j$ components of the $k$th column of $A$ by
the closure of $\alpha$ is assigned to the first $j$ components of $v$.\\
{\em Example 2.} $a(i,j)=a(i,j)\oplus a(i,1: n)\odot a(1: n,j)$ means that we add
($\oplus$) to the entry $a_{ij}$ of $A$ the result of the (universal) scalar multiplication of the $i$th row with the $j$th
column of $A$ (assumed that
$A$ is $n\times n$).\\
{\em Example 3.} $a(1: n,k)\odot b(l,1: m)$ means the outer product of
the $k$th column of $A$ with the $l$th row of $B$. The entries of resulting
matrix $C=(c_{ij})$ equal $c_{ij}=a_{ik}\odot b_{lj}$, for all $i=1,\ldots,n$
and $j=1,\ldots,m$.\\
{\em Example 4.} $x(1: n)\odot y(n: -1: 1)$ is the scalar product of
vector $x$ with vector $y$ whose components are taken in the reverse order:
the proper algebraic expression is $\bigoplus_{i=1}^n x_i\odot y_{n+1-i}$.\\
{\em Example 5.} The following cycle yields the same result as in the previous example:
$s=0$\\
{\bf for} $i=1:n$\\
$s=s\oplus x(i)\odot x(n+1-i)$\\
{\bf end}

\subsection{Escalator scheme and Gauss-Jordan elimination}
\label{ss:GJ}

We first analyse the basic escalator method, based on the definition
of matrix closures~\eqref{A_Star}. Let $A$ be a square matrix. Closures of its main submatrices $A_k$ can be found inductively, starting from
$A_1^*=(a_{11})^*$, the closure of the first diagonal entry. Generally we represent
$A_{k+1}$ as
\begin{equation*}
\label{ak+1}
A_{k+1}=
\begin{pmatrix}
A_k & g_k\\
h_k^T & a_{k+1}
\end{pmatrix},
\end{equation*}
assuming that we have found the closure of $A_k$. In this representation,
$g_k$ and $h_k$ are columns with $k$ entries and $a_{k+1}$ is a scalar. We also represent
$A_{k+1}^*$ as
\begin{equation*}
A_{k+1}^*=
\begin{pmatrix}
U_k & v_k\\
w_k^T & u_{k+1}
\end{pmatrix}.
\end{equation*}

Using~\eqref{A_Star} we obtain that

\begin{equation}
\label{bordering-mats}
\begin{split}
u_{k+1} &= (h_k^TA_k^*g_k\oplus a_{k+1})^*,\\
v_k &= A_k^* g_k u_{k+1},\\
w_k^T& = u_{k+1} h_k^T A_k^*,\\
U_k&= A_k^*g_ku_{k+1}h_k^TA_k^*\oplus A_k^*.
\end{split}
\end{equation}

An algorithm based on~\eqref{bordering-mats} can be written
as follows.

\begin{algorithm}
\label{a:bordering}
Escalator method for computing $A^*$
\end{algorithm}

\sssn {\bf Input:} an $n\times n$ matrix $A$ with entries $a(i,j)$,\\
also used to store the final result\\
and the intermediate results of the computation process.

\msn $a(1,1)=(a(1,1))^*$\\
{\bf for} $i=1: n-1$\\
$Ag=a(1: i,1: i)\odot a(1: i,i+1)$\\
$hA=a(i+1,1: i)\odot a(1: i,1: i)$\\
$a(i+1,i+1)=a(i+1,i+1)\oplus a(i+1,1: i)\odot Ag(1: i,1)$\\
$a(i+1,i+1)=(a(i+1,i+1))^*$\\
$a(1: i,i+1)=a(i+1,i+1)\odot Ag$\\
$a(i+1,1: i)=a(i+1,i+1)\odot hA$\\
$a(1: i,1: i) = a(1: i,1: i)\oplus Ag\odot a(i+1,i+1)\odot hA$\\
{\bf end}

\ssn In full analogy with its linear algebraic prototype,
the algorithm requires $n^3+O(n^2)$ operations of addition $\oplus$,
$n^3+O(n^2)$ operations of multiplication $\odot$, and $n$ operations of taking
algebraic closure. The linear-algebraic prototype of the method written above
is also called the {\em bordering method} in the
literature~\cite{Car-71,FF}.


Alternatively, we can obtain a solution of $X=AX\oplus B$ as a result of elimination
process, whose informal explanation is given below.
If $A^*$ is defined as $\bigoplus_{i\geq 0} A^i$ (including the scalar case), then $A^*B$ is the least
solution of $X=AX\oplus B$ for all $A$ and $B$ of appropriate sizes.
In this case, the solution found by the elimination process given below
coincides with $A^*B$.

For matrix $A=(a_{ij})$ and column vectors $x=(x_i)$ and $b=(b_i)$
(restricting without loss of generality to the column vectors),
the Bellman equation $x=Ax\oplus b$ can be written as

\begin{equation*}
\begin{pmatrix}
x_1\\
x_2\\
\vdots\\
x_n
\end{pmatrix}
=
\begin{pmatrix}
a_{11} & a_{12} &\ldots & a_{1n}\\
a_{21} & a_{22} &\ldots & a_{2n}\\
\vdots & \vdots &\ddots & \vdots\\
a_{n1} & a_{n2} & \ldots & a_{nn}
\end{pmatrix}
\begin{pmatrix}
x_1\\
x_2\\
\vdots\\
x_n
\end{pmatrix}
\oplus
\end{equation*}
\begin{equation}
\begin{pmatrix}
\1 & \0 & \ldots &\0 \\
\0 & \1 & \ldots &\0 \\
\vdots & \vdots &\ddots &\vdots\\
\0 &\0 &\ldots &\1
\end{pmatrix}
\begin{pmatrix}
b_1\\
b_2\\
\vdots\\
b_n
\end{pmatrix}.
\end{equation}

After expressing $x_1$ in terms of $x_2,\ldots,x_n$ from the first equation
and substituting this expression for $x_1$ in all other equations
from the second to the $n$th
we obtain

\begin{equation*}
\label{bellman-nice3}
\begin{pmatrix}
x_1\\
x_2\\
\vdots\\
x_n
\end{pmatrix}
=
\end{equation*}
\begin{equation*}
\begin{pmatrix}
\0 & (a_{11})^*a_{12} &\ldots & (a_{11})^*a_{1n}\\
\0 & a_{22}\oplus (a_{21}(a_{11})^* a_{12}) &\ldots & a_{2n}\oplus (a_{21}(a_{11})^* a_{1n})\\
\vdots & \vdots &\ddots & \vdots\\
\0 & a_{n2}\oplus (a_{n1}(a_{11})^* a_{12}) & \ldots &
a_{nn}\oplus (a_{n1}(a_{11})^* a_{1n})
\end{pmatrix}
\end{equation*}
\begin{equation}
\begin{pmatrix}
x_1\\
x_2\\
\vdots\\
x_n
\end{pmatrix}
\oplus
\begin{pmatrix}
(a_{11})^* & \0 & \ldots &\0 \\
a_{21}(a_{11})^*   & \1 & \ldots &\0 \\
\vdots & \vdots &\ddots &\vdots\\
a_{n1}(a_{11})^* &\0 &\ldots &\1
\end{pmatrix}
\begin{pmatrix}
b_1\\
b_2\\
\vdots\\
b_n
\end{pmatrix}
\end{equation}

Note that nontrivial entries in both matrices occupy complementary places,
so during computations both matrices can be stored in the same square array
$C^{(k)}$. Denote its elements by $c^{(k)}_{ij}$ where $k$ is the number of
eliminated variables. After $l-1$ eliminations we have

\begin{equation}
\label{e:xl}
\begin{split}
x_l&=(c_{ll}^{(l-1)})^*b_l,\\
c_{il}^{(l)}&=c_{il}^{(l-1)}(c_{ll}^{(l-1)})^*,\quad i=1,\ldots,l-1,l+1,\ldots,n\\
c_{ij}^{(l)}&=c_{ij}^{(l-1)}\oplus c_{il}^{(l-1)}(c_{ll}^{(l-1)})^* c_{lj}^{(l-1)},\\
i,j&=1,\ldots,l-1,l+1,\ldots,n\\
c_{li}^{(l)}&=(c_{ll}^{(l-1)})^*c_{li}^{(l-1)},\quad i=1,\ldots,l-1,l+1,\ldots,n\\
\end{split}
\end{equation}
After $n$ eliminations we get $x=C^{(n)}b$. Taking as $b$ any vector with one coordinate equal to $\1$ and
the rest equal to $\0$, we obtain $C^{(n)}=A^*$. We write out the following
algorithm based on recursion~\eqref{e:xl}.

\begin{algorithm}
\label{a:eliminate}
Gauss-Jordan elimination for computing $A^*$.
\end{algorithm}

\sssn {\bf Input}: an $n\times n$ matrix $A$ with entries $a(i,j)$,\\
also used to store the final result\\
and intermediate results of the computation process.

\msn {\bf for} $i=1: n$\\
$a(i,i)=(a(i,i))^*$\\
{\bf for} $k=1: n$\\
{\bf if} $k\neq i$\\
$a(k,i)=a(k,i)\odot a(i,i)$\\
{\bf end}\\
{\bf end}\\
{\bf for} $k=1: n$\\
{\bf for} $j=1: n$\\
{\bf if} $k\neq i\ \&\ j\neq i$\\
$a(k,j)=a(k,j)\oplus a(k,i)\odot a(i,j)$\\
{\bf end}\\
{\bf end}\\
{\bf for} $j=1: n$\\
{\bf if} $j\neq i$\\
$a(i,j)=a(i,i)\odot a(i,j)$\\
{\bf end}\\
{\bf end}\\
{\bf end}

\if{
\ssn TODO: 1) Quote a theorem stating that escalator and Gauss-Jordan
algorithms are equivalent or even identical;\\
2) Say that the Gauss-Jordan elimination is the same as the LDM decomposition
plus $M^*D^*L^*$ multiplication described below.
}\fi

\begin{remark}
{\rm Algorithm~\ref{a:eliminate} can be regarded as a ``universal
Floyd-Warshall algorithm'' generalizing the well-known algorithms
of Warshall and Floyd for computing the transitive closure of a
graph and all optimal paths on a graph. See, for instance,~\cite{Sedge} for
the description of these classical methods of discrete mathematics.
In turn, these methods can be regarded as specifications of
Algorithm~\ref{a:eliminate} to the cases of max-plus and Boolean semiring.}
\end{remark}

\begin{remark}
{\rm Algorithm~\ref{a:eliminate} is also
close to Yershov's ``refilling'' method for inverting matrices
and solving systems $Ax=b$ in the classical linear algebra,
see~\cite{FF}~Chapter~2 for details.}
\end{remark}

\subsection{Toeplitz systems}
\label{ss:Toep}

We start by considering the escalator method for finding the solution $x=A^*b$ to~
$x=Ax\oplus b$, where $x$ and $b$ are column
vectors. Firstly, we have $x^{(1)}=A_1^*b_1$. Let $x^{(k)}$ be the vector
found after $(k-1)$ steps, and let us write

\begin{equation*}
x^{(k+1)}=
\begin{pmatrix}
z\\
x_{k+1}
\end{pmatrix}.
\end{equation*}

Using~\eqref{bordering-mats} we obtain that
\begin{equation}
\label{bordering-vecs}
\begin{split}
x_{k+1}=u_{k+1}(h_k^T x^{(k)}\oplus b_{k+1}),\\
z=x^{(k)}\oplus A_k^*g_k x_{k+1}.
\end{split}
\end{equation}

We have to compute $A_k^* g_k$. In general, we would have to
use Algorithm~\ref{a:bordering}. Next we show that this calculation
can be done very efficiently when $A$ is symmetric Toep\-litz.

Formally, a matrix $A\in\Mat_{nn}(\semr)$ is called {\em Toep\-litz} if there exist
scalars $r_{-n+1},\ldots,r_0,\ldots,r_{n-1}$ such that $A_{ij}=r_{j-i}$ for all $i$ and $j$.
Informally, Toep\-litz matrices are such that their entries are constant along any line parallel to the
main diagonal (and along the main diagonal itself). For example,

\begin{equation*}
A=
\begin{pmatrix}
r_0 & r_1 & r_2 & r_3\\
r_{-1} & r_0 & r_1 & r_2\\
r_{-2} & r_{-1} & r_0 & r_1\\
r_{-3} & r_{-2} & r_{-1} & r_0
\end{pmatrix}
\end{equation*}
is Toeplitz. Such matrices are not necessarily symmetric. However, they are always
{\em persymmetric}, that is, symmetric with respect to the inverse diagonal. This property
is algebraically expressed as $A=E_nA^TE_n$, where $E_n=[e_n,\ldots,e_1]$.
By $e_i$ we denote the column whose $i$th entry is $\bunity$ and other entries are
$\bzero$. The property $E_n^2=I_n$ (where $I_n$ is the $n\times n$ identity
matrix) implies that the product of two persymmetric matrices is persymmetric. Hence
any degree of a persymmetric matrix is persymmetric, and so is the closure of a
persymmetric matrix. Thus, if $A$ is persymmetric, then
\begin{equation}
\label{pers}
E_nA^*=(A^*)^TE_n.
\end{equation}

Further we deal only with symmetric Toeplitz matrices. Consider the equation
$y=T_ny\oplus r^{(n)}$, where $r^{(n)}=(r_1,\ldots,r_n)^T$, and $T_n$ is defined by the scalars
$r_0,r_1,\ldots,r_{n-1}$ so that $T_{ij}=r_{|j-i|}$ for all $i$ and $j$. This is a generalization of the
Yule-Walker problem~\cite{GvL}. Assume that we have obtained the least solution $y^{(k)}$
to the system $y=T_k y\oplus r^{(k)}$ for some $k$ such that $1\leq k\leq n-1$, where $T_k$ is the main
$k\times k$ submatrix of $T_n$. We write $T_{k+1}$ as
$$ T_{k+1}=\begin{pmatrix}
T_k & E_k r^{(k)}\\
r^{(k)T} E_k & r_0
\end{pmatrix}. $$
We also write $y^{(k+1)}$ and $r^{(k+1)}$ as

\begin{equation*}
y^{(k+1)}=\begin{pmatrix}
z\\
\alpha_k
\end{pmatrix},\quad
r^{(k+1)}=
\begin{pmatrix}
r^{(k)}\\
r_{k+1}
\end{pmatrix}.
\end{equation*}

Using~\eqref{bordering-vecs}, \eqref{pers} and the identity $T_k^*r^{(k)}=y^{(k)}$,
we obtain that

\begin{equation*}
\begin{split}
\alpha_k &= (r_0\oplus r^{(k)T}y^{(k)})^*(r^{(k)T}E_ky^{(k)}\oplus r_{k+1}),\\
z &= E_k y^{(k)}\alpha_k\oplus y^{(k)}.
\end{split}
\end{equation*}

Denote $\beta_k=r_0\oplus r^{(k)T}y^{(k)}$. The following argument shows that $\beta_k$
can be found recursively if $(\beta_{k-1}^*)^{-1}$ exists.

\begin{equation*}
\begin{split}
\beta_k &= r_0\oplus [r^{(k-1)T}\ r_k]
\begin{pmatrix}
E_{k-1}y^{(k-1)}\alpha_{k-1}\oplus y^{(k-1)}\\
\alpha_{k-1}
\end{pmatrix}\\
&=r_0\oplus r^{(k-1)T}y^{(k-1)}\oplus
(r^{(k-1)T} E_{k-1} y^{(k-1)}\oplus
\end{split}
\end{equation*}
\begin{equation}
\label{betak}
r_k)\alpha_{k-1} = \beta_{k-1}\oplus(\beta_{k-1}^*)^{-1}\odot\alpha_{k-1}^2.
\end{equation}

Existence of $(\beta_{k-1}^*)^{-1}$ is not universal, and this will make us write two versions of our algorithm,
the first one involving~\eqref{betak}, and the second one not involving it. We will write these
two versions in one program and mark the expressions which refer only to the first version
or to the second one by the MATLAB-style comments $\% 1$ and $\% 2$, respectively.
Collecting the expressions for $\beta_k$, $\alpha_k$ and $z$ we obtain the following
recursive expression for $y^{(k)}$:

\begin{equation}
\label{durbin}
\begin{split}
\beta_k &= r_0\oplus r^{(k)T} y^{(k)},\qquad \% 2\\
\beta_k &=\beta_{k-1}\oplus(\beta_{k-1}^*)^{-1}\odot\alpha_{k-1}^2,\qquad \% 1\\
\alpha_k &= (\beta_k)^*\odot(r^{(k)T}E_ky^{(k)}\oplus r_{k+1}),\\
y^{(k+1)}&=
\begin{pmatrix}
E_k y^{(k)}\alpha_k\oplus y^{(k)}\\
\alpha_k
\end{pmatrix}.
\end{split}
\end{equation}
Recursive expression~\eqref{durbin} is a generalized version of the Durbin method for the
Yule-Walker problem, see~\cite{GvL} Algorithm 4.7.1 for a prototype.

\begin{algorithm}
\label{a:YW}
The Yule-Walker problem for the Bellman equations with symmetric
Toeplitz matrix.
\end{algorithm}

\sssn {\bf Input:} $r_0$: scalar,\\
$r$: $n-1\times 1$ vector;

\msn $y(1)=r_0^*\odot r(1)$\\
$\beta=r_0\qquad \% 1$\\
$\alpha=r_0^*\odot r(1)$\\
{\bf for} $k=1: n-1$\\
$\beta=r_0\oplus r(1: k)\odot y(1: k)\qquad \% 2$\\
$\beta=\beta\oplus(\beta^*)^{-1}\odot\alpha^2\qquad \% 1$\\
$\alpha=\beta^*\odot(r(k: -1: 1)\odot y(1: k)\oplus r(k+1))$\\
$z(1: k)=y(1: k)\oplus y(k: -1: 1)\odot\alpha$\\
$y(1: k)=z(1: k)$\\
$y(k+1)=\alpha$\\
{\bf end}

\msn {\bf Output:} vector $y$.

\ssn In the general case, the algorithm requires $3/2 n^2+O(n)$ operations $\oplus$ and
$\odot$ each, and just $n^2+O(n)$
of $\oplus$ and $\odot$ if inversions of algebraic closures are
allowed (as usual, just $n$ such closures are required in both cases).

Now we consider the problem of finding $x^{(n)}=T_n^*b^{(n)}$ where $T_n$ is as above and
$b^{(n)}=(b_1,\ldots,b_n)$ is arbitrary. We also introduce the column vectors $y^{(k)}$ which solve the Yule-Walker problem:
$y^{(k)}=T_k^*r^{(k)}$. The main idea is to find the expression for $x^{(k+1)}=T_{k+1}^*b^{(k+1)}$ involving
$x^{(k)}$ and $y^{(k)}$. We write $x^{(k+1)}$ and $b^{(k+1)}$ as
\begin{equation*}
x^{(k+1)}=
\begin{pmatrix}
v\\
\mu_k
\end{pmatrix},\quad
b^{(k+1)}=
\begin{pmatrix}
b^{(k)}\\
b_{k+1}
\end{pmatrix}.
\end{equation*}

Making use of the persymmetry of $T_k^*$ and of the identities $T_k^*b_k=x^{(k)}$ and $T_k^*r_k=y^{(k)}$,
we specialize expressions~\eqref{bordering-vecs} and obtain that

\begin{equation*}
\begin{split}
\mu_k &= (r_0\oplus r^{(k)T} y^{(k)})^*\odot (r^{(k)T}E_kx^{(k)}\oplus b_{k+1}),\\
v &= E_k y^{(k)}\mu_k\oplus x^{(k)}.
\end{split}
\end{equation*}

The coefficient $r_0\oplus r^{(k)T}y^{(k)}=\beta_k$ can be expressed again as
$\beta_k=\beta_{k-1}\oplus(\beta_{k-1}^*)^{-1}\odot(\alpha_{k-1})^2$, if the closure $(\beta_{k-1})^*$ is invertible.
Using this we obtain the following recursive expression:

\begin{equation}
\label{levinson}
\begin{split}
\beta_k &= r_0\oplus r^{(k)T} y^{(k)},\qquad \% 2\\
\beta_k &=\beta_{k-1}\oplus(\beta_{k-1}^*)^{-1}\odot\alpha_{k-1}^2,\qquad \% 1\\
\mu_k &= \beta_k^*\odot(r^{(k)T}E_kx^{(k)}\oplus b_{k+1}),\\
x^{(k+1)} &=
\begin{pmatrix}
E_k y^{(k)}\mu_k\oplus x^{(k)}\\
\mu_k
\end{pmatrix}.
\end{split}
\end{equation}

Expressions~\eqref{durbin} and~\eqref{levinson} yield the following generalized version of the Levinson algorithm
for solving linear symmetric Toeplitz systems,
see~\cite{GvL} Algorithm~4.7.2 for a prototype:

\begin{algorithm}
\label{a:levinson}
Bellman system with symmetric Toeplitz matrix
\end{algorithm}

\sssn {\bf Input:} $r_0$: scalar,\\
$r$: $1\times n-1$ row vector;\\
$b$: $n\times 1$ column vector.

\msn $y(1)=r_0^*\odot r(1)$;\quad $x(1)=r_0^*\odot b(1)$;\\
$\beta=r_0\qquad \% 1$\\
$\alpha=r_0^*\odot r(1)$\\
{\bf for} $k=1: n-1$\\
$\beta=r_0\oplus r(1: k)\odot y(1: k)\qquad \% 2$\\
$\beta=\beta\oplus(\beta^*)^{-1}\odot\alpha^2\qquad \% 1$\\
$\mu=\beta^*\odot(r(k: -1: 1)\odot x(1: k)\oplus b(k+1))$\\
$v(1: k)=x(1: k)\oplus y(k: -1: 1)\odot\mu$\\
$x(1: k)=v(1: k)$\\
$x(k+1)=\mu$\\
{\bf if} $k<n-1$\\
$\alpha=\beta^*\odot(r(k: -1: 1)\odot y(1: k)\oplus r(k+1))$\\
$z(1: k)=y(1: k)\oplus y(k: -1: 1)\odot\alpha$\\
$y(1: k)=z(1: k)$\\
$y(k+1)=\alpha$\\
{\bf end}\\
{\bf end}

\msn {\bf Output:} vector $x$.

\ssn In the general case, the algorithm requires $5/2 n^2+O(n)$ operations $\oplus$ and
$\odot$ each, and just $2n^2+O(n)$
of $\oplus$ and $\odot$ if inversions of algebraic closures are
allowed (as usual, just $n$ such closures are required in both cases).

\subsection{LDM decomposition}
\label{ss:LDM}

Factorization of a matrix into the product $A = LDM$, where $L$ and $M$
are lower and upper triangular matrices with a unit diagonal,
respectively, and $D$ is a diagonal matrix, is used for solving
matrix equations $AX = B$. We construct a similar
decomposition for the Bellman equation $X = AX \oplus B$.

For the case $AX = B$, the decomposition $A = LDM$ induces the following
decomposition of the initial equation:
\begin{equation}
   LZ = B, \qquad DY = Z, \qquad MX = Y.
\end{equation}
Hence, we have
\begin{equation}
   A^{-1} = M^{-1}D^{-1}L^{-1},
\label{AULinv}
\end{equation}

if $A$ is invertible. In essence, it is sufficient to find the
matrices $L$, $D$ and $M$, since the linear system $AX=B$ is easily
solved by a combination of the forward substitution for $Z$, the
trivial inversion of a diagonal matrix for $Y$, and the back
substitution for $X$.

Using the LDM-factorization of $AX=B$ as a prototype, we can write

\begin{equation}
   Z = LZ \oplus B, \qquad Y = DY \oplus Z, \qquad X = MX \oplus Y.
\label{LDM}
\end{equation}
Then
\begin{equation}
   A^* = M^*D^*L^*.
\label{AMDLstar}
\end{equation}

A triple $(L,D,M)$ consisting of a lower triangular, diagonal, and
upper triangular matrices is called an $LDM$-{\it factorization} of a
matrix $A$ if relations~\eqref{LDM} and~\eqref{AMDLstar}
are satisfied. We note that
in this case, the principal diagonals of $L$ and $M$ are zero.

Our universal modification of the $LDM$-fac\-tor\-iza\-tion used in matrix
analysis for the equation $AX=B$ is similar to the $LU$-factorization of Bellman equation suggested by Carr\'e in~\cite{Car-71,Car:79}.

\if{
We stress that the algorithm described in this section
can be applied to matrix
computations over any semiring under the condition that the unary
operation $a\mapsto a^*$ is applicable every time it is encountered
in the computational process. Indeed, when constructing the
algorithm, we use only the basic semiring operations of addition
$\oplus$ and multiplication $\odot$ and the properties of
associativity, commutativity of addition, and distributivity of
multiplication over addition.
}\fi

If $A$ is a symmetric matrix over a semiring with a commutative
multiplication, the amount of computations can be halved, since
$M$ and $L$ are mapped into each other under transposition.

We begin with the case of a triangular matrix $A = L$ (or $A = M$).
Then, finding $X$ is reduced to the forward (or back) substitution.
Note that in this case, equation $X=AX\oplus B$ has unique solution, which
can be found by the obvious algorithms given below. In these algorithms $B$ is a vector (denoted by $b$), however they could be modified to the case when
$B$ is a matrix of any appropriate size. We are interested only in the case
of strictly lower-triangular, resp. strictly upper-triangular matrices, when
$a_{ij}=0$ for $i\leq j$, resp. $a_{ij}=0$ for $i\geq j$.

\begin{algorithm}
\label{lower-triang}
Forward substitution.
\end{algorithm}

\sssn {\bf Input:} Strictly lower-triangular $n\times n$ matrix $l$;\\
$n\times 1$ vector $b$.

\msn {\bf for} $k=2: n$\\
$y(k)=l(k,1: k-1)\odot y(1:k-1)$\\
{\bf end}

\msn {\bf Output:} vector $y$.

\begin{algorithm}
\label{upper-triang}
Backward substitution.
\end{algorithm}

\sssn {\em Input:} Strictly upper-triangular $n\times n$ matrix $m$;\\
$n\times 1$ vector $b$.

\msn {\bf for} $k=n-1: -1: 1$\\
$y(k)=m(k,k+1: n)\odot y(k+1: n)$\\
{\bf end}

\msn {\bf Output:} vector $y$.

\ssn Both algorithms require $n^2/2+O(n)$ operations $\oplus$ and $\odot$,
and no algebraic closures.

After performing a LDM-decomposition we also need to compute
the closure of a diagonal matrix: this is done entrywise.

We now proceed with the algorithm of LDM decomposition itself, that is,
computing matrices $L$, $D$ and $M$ satisfying~\eqref{LDM}
and~\eqref{AMDLstar}. First we give an algorithm, and then we proceed
with its explanation.

\begin{algorithm}
\label{a:ldm}
LDM-decomposition (version 1).
\end{algorithm}

\sssn {\bf Input}: an $n\times n$ matrix $A$ with entries $a(i,j)$,\\
also used to store the final result\\
and intermediate results of the computation process.

\msn {\bf for} $j=1: n-1$\\
$v(j)=(a(j,j))^*$\\
$a(j+1: n,j)=a(j+1: n,j)\odot v(j)$\\
$a(j+1: n,j+1: n)=a(j+1: n,j+1: n)\oplus a(j+1: n,j)\odot a(j,j+1: n)$\\
$a(j,j+1: n)=v(j)\odot a(j,j+1: n)$\\
{\bf end}

\ssn The algorithm requires $n^3/3+O(n^2)$ operations $\oplus$ and $\odot$,
and $n-1$ operations of algebraic closure.

The strictly triangular matrix $L$ is written in the lower triangle,
the strictly upper triangular matrix $M$ in the upper triangle,
and the diagonal matrix $D$ on the diagonal of the matrix computed by
Algorithm~\ref{a:ldm}. We now show that $A^*=M^*D^*L^*$. Our argument
is close to that of~\cite{BC-75}.

We begin by representing, in analogy with the escalator method,

\begin{equation}
\label{a-repr}
A=
\begin{pmatrix}
a_{11} & h^{(1)}\\
g^{(1)} & B^{(1)}
\end{pmatrix}
\end{equation}

It can be verified that

\begin{equation*}
A^*=
\begin{pmatrix}
\bunity & h^{(1)} a_{11}^*\\
\bzero_{n-1\times 1} & I_{n-1}
\end{pmatrix}\odot
\end{equation*}
\begin{equation}
\label{a*ldm1}
\begin{pmatrix}
a_{11}^* & \bzero_{1\times n-1} \\
\bzero_{n-1\times 1} & (h^{(1)}a_{11}^*g^{(1)}\oplus B^{(1)})^*
\end{pmatrix}
\begin{pmatrix}
\bunity & \bzero_{1\times n-1}\\
a_{11}^*g^{(1)} & I_{n-1}
\end{pmatrix}
\end{equation}
as the multiplication on the right hand side leads to expressions fully analogous
to~\eqref{bordering-mats}, where\\
$(h^{(1)}a_{11}^*g^{(1)}\oplus B^{(1)})^*$ plays the role of
$u_{k+1}$. Here and in the sequel, $O_{k\times l}$ denotes the $k\times l$ matrix consisting only of zeros,
and $I_l$ denotes the identity matrix of size $l$.
This can be also rewritten as

\begin{equation}
\label{a*ldm2}
A^*=M_1^*D_1^*(A^{(2)})^* L_1^*,
\end{equation}

where

\begin{equation*}
M_1=
\begin{pmatrix}
\bzero & h^{(1)} a_{11}^*\\
\bzero_{(n-1)\times 1} & \bzero_{(n-1)\times(n-1)}
\end{pmatrix},
\end{equation*}
\begin{equation*}
D_1=
\begin{pmatrix}
a_{11} & \bzero_{1\times (n-1)}\\
\bzero_{(n-1)\times 1} & \bzero_{(n-1)\times(n-1)}
\end{pmatrix},
\end{equation*}
\begin{equation*}
A^{(2)} =
\begin{pmatrix}
\bzero_{1\times 1} & \bzero_{1\times(n-1)}\\
\bzero_{(n-1)\times 1} & R^{(2)}
\end{pmatrix},
\end{equation*}
\begin{equation*}
\label{l1d1m1}
L_1 =
\begin{pmatrix}
\bzero_{1\times 1} & \bzero_{1\times (n-1)}\\
a_{11}^*g^{(1)} & O_{(n-1)\times (n-1)}
\end{pmatrix},
\end{equation*}

\begin{equation}
R^{(2)}=h^{(1)}a_{11}^*g^{(1)}\oplus B^{(1)}.
\end{equation}

Here we used in particular that $L_1^2=0$ and $M_1^2=0$ and hence
$L_1^*=I\oplus L_1$ and $M_1^*=I\oplus M_1$.

The first step of Algorithm~\ref{a:ldm} ($k=1$) computes

\begin{equation}
\label{s1}
\begin{pmatrix}
a_{11} & h^{(1)} a_{11}^*\\
a_{11}^*g^{(1)} & R^{(2)}
\end{pmatrix}=A^{(2)}\oplus L_1\oplus M_1\oplus D_1,
\end{equation}

which contains all relevant information.

We can now continue with the submatrix $R^{(2)}$ of $A^{(2)}$ factorizing it
as in~\eqref{a*ldm1} and~\eqref{a*ldm2}, and so on. Let us now formally
describe the
$k$th step of this construction, corresponding to the $k$th step of
Algorithm~\ref{a:ldm}. On that general step we deal with

\begin{equation}
\label{e:ak}
A^{(k)}=
\begin{pmatrix}
\bzero_{(k-1)\times(k-1)} & \bzero_{(k-1)\times(n-k+1)}\\
\bzero_{(n-k+1)\times(k-1)} & R^{(k)}
\end{pmatrix},
\end{equation}

where

\begin{equation*}
\label{e:rk}
R^{(k)}= h^{(k-1)}(a_{k-1,k-1}^{(k-1)})^*g^{(k-1)}\oplus B^{(k-1)}=
\end{equation*}
\begin{equation}
\begin{pmatrix}
a_{kk}^{(k)} & h^{(k)}\\
g^{(k)} & B^{(k)}
\end{pmatrix}.
\end{equation}

Like on the first step we represent

\begin{equation}
\label{a*ldmk1}
(A^{(k)})^*=M_k^*D_k^*(A^{(k+1)})^*L_k^*,
\end{equation}

where

\begin{equation}
\begin{split}
\label{a*ldmk2}
M_k&=
\begin{pmatrix}
\bzero_{(k-1)\times(k-1)} & \bzero_{(k-1)\times 1} & \bzero_{(k-1)\times(n-k)}\\
\bzero_{1\times (k-1)} & \bzero_{1\times 1} & h^{(k)}(a_{kk}^{(k)})^*\\
\bzero_{(n-k)\times(k-1)} & \bzero_{(n-k)\times 1} & \bzero_{(n-k)\times(n-k)}
\end{pmatrix},\\
D_k&=
\begin{pmatrix}
\bzero_{(k-1)\times (k-1)} & \bzero_{(k-1)\times 1} & \bzero_{(k-1)\times(n-k)}\\
\bzero_{1\times(k-1)} & a_{kk}^{(k)} & \bzero_{1\times(n-k)}\\
\bzero_{(n-k)\times (k-1)} & \bzero_{(n-k)\times 1} & \bzero_{(n-k)\times (n-k)}
\end{pmatrix},\\
L_k&=
\begin{pmatrix}
\bzero_{(k-1)\times (k-1)} & \bzero_{(k-1)\times 1} & \bzero_{(k-1)\times(n-k)}\\
\bzero_{1\times(k-1)} & \bzero_{1\times 1} & \bzero_{1\times(n-k)}\\
\bzero_{(n-k)\times (k-1)} & (a_{kk}^{(k)})^*g^{(k)} & \bzero_{(n-k)\times (n-k)}
\end{pmatrix},\\
A^{(k+1)}&=
\begin{pmatrix}
\bzero_{k\times k} & \bzero_{k\times(n-k)}\\
\bzero_{(n-k)\times k} & R^{(k+1)}
\end{pmatrix},\\
R^{(k+1)}&=h^{(k)}(a_{kk}^{(k)})^*g^{(k)}\oplus B^{(k)}.
\end{split}
\end{equation}
Note that we have the following recursion for the entries of
$A^{(k)}$:
\begin{equation}
\label{akrec}
a_{ij}^{(k+1)}=
\begin{cases}
\0, & \text{if $i\leq k$ or $j\leq k$},\\
a_{ij}^{(k)}\oplus a_{ik}^{(k)}(a_{kk}^{(k)})^*a_{kj}^{(k)}, & \text{otherwise}.
\end{cases}
\end{equation}
This recursion is immediately seen in Algorithm~\ref{a:ldm}.
Moreover it can be shown by induction that the matrix computed on
the $k$th step of that algorithm equals
\begin{equation}
A^{(k+1)}\oplus\bigoplus_{i=1}^k L_i \oplus\bigoplus_{i=1}^k M_i\oplus
\bigoplus_{i=1}^k D_i.
\end{equation}
In other words, this matrix is composed from $h^{(1)}a_{11}^*$, ...,
$h^{(k)}(a_{kk}^{(k)})^*$ (in the upper triangle),
$a_{11}^*g^{(1)}$, ..., $(a_{kk}^{(k)})^*g^{(k)}$ (in the lower triangle),
$a_{11},\ldots,a_{kk}^{(k)}$ (on the diagonal), and $R^{(k+1)}$ (in the south-eastern corner).

After assembling and unfolding all expressions \eqref{a*ldmk1} for $A^{(k)}$,
where $k=1,\ldots,n$, we obtain
\begin{equation}
\label{e:ldm1}
A^*=M_1^*D_1^*\cdots M_n^*D_n^*L_n^*\cdots L_1^*.
\end{equation}
(actually, $M_n=L_n=0$ and hence $M_n^*=L_n^*=I$).
Noticing that $D_i^*$ and $M_j^*$ commute for $i<j$ we can rewrite
\begin{equation}
\label{e:ldm2}
A^*=M_1^*\cdots M_n^* D_1^*\cdots D_n^* L_n^*\cdots L_1^*.
\end{equation}
Consider the identities
\begin{equation}
\label{sumprod}
\begin{split}
(D_1\oplus\ldots\oplus D_n)^*&=D_1^*\cdots D_n^*,\\
(L_1\oplus\ldots\oplus L_n)^*&=L_n^*\cdots L_1^*,\\
(M_1\oplus\ldots\oplus M_n)^*&=M_1^*\cdots M_n^*.
\end{split}
\end{equation}
The first of these identities is evident. For the other two, observe
that $M_k^2=L_k^2=0$ for all $k$, hence $M_k^*=I\oplus M_k$ and
$L_k^*=I\oplus L_k$. Further, $L_iL_j=0$ for $i>j$ and
$M_iM_j=0$ for $i<j$. Using these identities it can be shown that
\begin{equation}
\label{sumprod-exp}
\begin{split}
&(L_1\oplus\ldots\oplus L_n)^*=\bigoplus_{i=0}^{n-1} (L_1\oplus\ldots\oplus L_n)^i=\\
&=(I\oplus L_n)\cdots (I\oplus L_1)=L_n^*\cdots L_1^*,\\
&(M_1\oplus\ldots\oplus M_n)^*=\bigoplus_{i=0}^{n-1} (M_1\oplus\ldots\oplus M_n)^i=\\
&=(I\oplus M_1)\cdots (I\oplus M_n)=M_1^*\cdots M_n^*,
\end{split}
\end{equation}
which yields the last two identities of~\eqref{sumprod}. Notice that in~\eqref{sumprod-exp} we have used the nilpotency of $L_1\oplus\ldots\oplus L_n$ and
$M_1\oplus\ldots \oplus M_n$, which allows to apply~\eqref{a*nilp}.

It can be seen that the matrices
$M:=M_1\oplus\ldots\oplus M_n$, $L:=L_1\oplus\ldots\oplus L_n$ and
$D:=D_1\oplus\ldots\oplus D_n$ are contained in the upper triangle,
in the lower triangle and, respectively, on the diagonal of the matrix
computed by Algorithm~\ref{a:ldm}. These matrices satisfy the LDM
decomposition $A^*=M^*D^*L^*$. This concludes the explanation of
Algorithm~\ref{a:ldm}.

In terms of matrix computations, Algorithm~\ref{a:ldm} is a
version of LDM decomposition with outer product.
This algorithm can be reorganized to make it almost
identical with~\cite{GvL}, Algorithm 4.1.1:

\begin{algorithm}
\label{a:ldmgvl}
LDM-decomposition (version 2).
\end{algorithm}

\sssn {\bf Input:} an $n\times n$ matrix $A$ with entries $a(i,j)$,\\
also used to store the final result\\
and intermediate results of the computation process.

\msn {\bf for} $j=1: n$\\
$v(1:j)=a(1:j,j)$\\
{\bf for} $k=1: j-1$
$v(k+1: j)=v(k+1:j)\oplus a(k+1: j,k)\odot v(k)$\\
{\bf end}\\
{\bf for} $i=1: j-1$\\
$a(i,j)=(a(i,i))^*\odot v(i)$\\
{\bf end}\\
$a(j,j)=v(j)$\\
{\bf for} $k=1: j-1$\\
$a(j+1: n,j) = a(j+1: n,j)\oplus a(j+1: n,k)\odot v(k)$\\
{\bf end}\\
$d=(v(j))^*$\\
$a(j+1: n,j)=a(j+1: n,j)\odot d$\\
{\bf end}

\ssn This algorithm performs exactly the same operations as Algorithm~\ref{a:ldm},
computing consecutively one column of the result after another. Namely,
in the first half of the main loop it computes the entries $a_{ij}^{(i)}$ for $i=1,\ldots,j$, first under the guise of the entries of $v$
and finally in the assignment ``$a(i,j)=(a(i,i))^*\odot v(i)$''.  In the
second half of the main loop it computes $a_{kj}^{(j)}$. The complexity of this
algorithm is the same as that of Algorithm~\ref{a:ldm}.

\subsection{LDM decomposition with symmetry and band structure}
\label{ss:LDMspec}

When matrix $A$ is symmetric, that is, $a_{ij}=a_{ji}$ for all $i,j$,
it is natural to expect that LDM decomposition must be symmetric too, that is,
$M=L^T$. Indeed, going through the reasoning of the previous section, it can be
shown by induction that all intermediate matrices $A^{(k)}$ are symmetric,
hence $M_k=L_k^T$ for all $k$ and $M=L^T$. We now present two versions of
symmetric LDM decomposition, corresponding to the two versions of LDM decomposition
given in the previous section. Notice that the amount of computations in these
algorithms is nearly halved with respect to their full versions. In both cases
they require $n^3/6+O(n^2)$ operations $\oplus$ and $\odot$(each) and $n-1$
operations of taking algebraic closure.

\begin{algorithm}
\label{a:ldl}
Symmetric LDM-decomposition (version 1).
\end{algorithm}

\sssn {\bf Input:} an $n\times n$ symmetric matrix $A$ with entries $a(i,j)$,\\
also used to store the final result\\
and intermediate results of the computation process.

\msn {\bf for} $j=1: n-1$\\
$v(j)=(a(j,j))^*$\\
{\bf for} $k=j+1: n$\\
{\bf for} $l=j+1: k$\\
$a(k,l) = a(k,l)\oplus a(k,j)\odot v(j)\odot a(l,j)$\\
{\bf end}\\
{\bf end}\\
$a(j+1: n,j)= a(j+1: n,j)\odot v(j)$\\
{\bf end}

\ssn The strictly triangular matrix $L$ is contained in the lower triangle
of the result, and the matrix $D$ is on the diagonal.

The next version generalizes~\cite{GvL} Algorithm 4.1.2.
Like in the prototype, the idea is to use the symmetry of $A$
precomputing the first $j-1$ entries of $v$ inverting the assignment
``$a(i,j)=a(i,i)^*\odot v(i)$'' for $i=1,\ldots, j-1$. This is possible
since $a(j,i)=a(i,j)$ belong to the first $j-1$ columns
of the result that have been computed on the previous stages.

\begin{algorithm}
\label{a:ldlgvl}
Symmetric LDM-decomposition\\ (version 2).
\end{algorithm}

\sssn $A$ is an $n\times n$ symmetric matrix with entries $a(i,j)$,\\
also used to store the final result\\
and intermediate results of the computation process.

\msn {\bf for} $j=1: n$\\
{\bf for} $i=1: j-1$\\
$v(i)=(a(i,i)^*)^{-1}a(i,j)$\\
{\bf end}\\
$v(j)=v(j)\oplus a(j,1: j-1)\odot v(1: j-1)$\\
$a(j,j)=v(j)$\\
{\bf for} $k=1: j-1$\\
$a(j+1: n,j) = a(j+1: n,j)\oplus a(j+1: n,k)\odot v(k)$\\
{\bf end}\\
$d=(v(j))^*$\\
$a(j+1: n,j)=a(j+1: n,j)\odot d$
{\bf end}

\ssn Note that this version requires invertibility of the closures
$a(i,i)^*$ computed by the algorithm.

\begin{remark}
In the case of idempotent semiring we have $(D^*)^2=D^*$, hence\\
$A^*=(M^*D^*)(D^*L^*)$. When $A$ is symmetric we can write $A^*=(G^*)^TG^*$
where $G=D^*L$. Evidently, this  {\bf idempotent Cholesky factorization} can be
computed by minor modifications of Algorithms~\ref{a:ldl} and~\ref{a:ldlgvl}.
See also~\cite{GvL}, Algorithm 4.2.2.
\end{remark}

$A=(a_{ij})$ is called a {\bf band matrix} with upper bandwidth $q$ and lower
bandwidth $p$ if $a_{ij}=0$ for all $j>i+q$ and all $i>j+p$. A band matrix
with $p=q=1$ is called {\bf tridiagonal}. To generalize a specific LDM decomposition
with band matrices, we need to show that the band parameters of the matrices $A^{(2)},\ldots, A^{(n)}$ computed in the process of LDM decomposition are not
greater than the parameters of $A^{(1)}=A$. Assume by induction that
$A=A^{(1)},\ldots,A^{(k)}$ have the required band parameters, and consider an entry
$a_{ij}^{(k+1)}$ for $i>j+p$. If $i\leq k$ or $j\leq k$ then $a_{ij}^{(k+1)}=\0$,
so we can assume $i>k$ and $j>k$. In this case $i>k+p$, hence $a_{ik}^{(k)}=\0$
and
$$a_{ij}^{(k+1)}=a_{ij}^{(k)}\oplus a_{ik}^{(k)}(a_{kk}^{(k)})^*a_{kj}^{(k)}=\0.$$
Thus we have shown that the lower bandwidth of $A^{(k)}$ is not greater than $p$.
It can be shown analogously that its upper bandwidth does not exceed $q$.
We use this to construct the following band version of LDM decomposition,
see \cite{GvL} Algorithm 4.3.1 for a prototype.

\begin{algorithm}
\label{a:ldmband}
LDM decomposition of a band matrix.
\end{algorithm}

\sssn $A$ is an $n\times n$ band matrix with entries $a(i,j)$,\\
lower bandwidth $p$ and upper bandwidth $q$\\
also used to store the final result\\
and intermediate results of the computation process.

\msn {\bf for} $j=1: n-1$\\
$v(j)=(a(j,j))^*$\\
{\bf for} $i=j+1: \min(j+p,n)$\\
$a(i,j)=a(i,j)\odot v(j)$\\
{\bf end}\\
{\bf for} $k=j+1: \min(j+q,n)$\\
{\bf for} $i=j+1: \min(j+p,n)$\\
$a(k,j)=a(k,j)\oplus a(k,i)\odot a(i,j)$\\
{\bf end}\\
{\bf end}\\
{\bf for} $k=j+1: \min(j+q,n)$\\
$a(j,k)=v(j)\odot a(j,k)$\\
{\bf end}\\
{\bf end}

\ssn When $p$ and $q$ are fixed and $n>>p,q$ is variable,
it can be seen that the algorithm performs
approximately $npq$ operations $\odot$ and $\oplus$ each.

\begin{remark}
{\em There are important special kinds of band matrices, for instance, Hessenberg
and tridiagonal matrices. Hessenberg matrices are defined as band matrices with
$p=1$ and $q=n$, while in the case of tridiagonal matrices $p=q=1$. It is
straightforward to write further adaptations of Algorithm~\ref{a:ldmband}
to these cases.}
\end{remark}

\subsection{Iteration schemes}
\label{ss:Iter}

We are not aware of any truly universal scheme, since the decision when
such schemes work and when they should be stopped depends both
on the semiring and on the representation of data.

Our first scheme is derived from the following iteration
process:
\begin{equation}
\label{e:jacobi}
X^{(k+1)}=AX^{(k)}\oplus B
\end{equation}
trying to solve the Bellman equation $X=AX\oplus B$. Iterating expressions
~\eqref{e:jacobi} for all $k$ up to $m$ we obtain
\begin{equation}
\label{e:jac-unfold}
X^{(m)}=A^m X^{(0)}\oplus \bigoplus_{i=0}^{m-1} A^iB
\end{equation}
Thus the result crucially depends on the behaviour of $A^mX^{(0)}$.
The algorithm can be written as follows (for the case when $B$ is a column vector).

\begin{algorithm}
\label{a:jacobi}
Jacobi iterations
\end{algorithm}

\sssn {\bf Input:} $n\times n$ matrix $A$ with entries $a(i,j)$;\\
$n\times 1$ column vectors $b$ and $x$

\msn situation$=$'proceed'\\
{\bf while}  situation$==$'proceed'\\
$x= A\odot x\oplus b$\\
situation$=${\em newsituation}(...)\\
{\bf if} situation$==$'no convergence'\\
{\bf disp}('Jacobi iterations did not converge')\\
{\em exit}\\
{\bf end}\\
{\bf if} situation$==$'convergence'\\
{\bf disp}('Jacobi iterations converged')\\
{\em exit}\\
{\bf end}\\
{\bf end}

\msn {\bf Output:} situation, $x$.

\ssn Next we briefly discuss the behaviour of Jacobi iteration scheme
over the usual arithmetic with nonnegative real numbers, and over
semiring $\rmax$. For simplicity, in both cases we restrict to the case
of {\em irreducible} matrix $A$, that is, when the associated digraph is strongly
connected.

Over the usual arithmetic, it is well known that (in the irreducible nonnegative case) the Jacobi iterations converge if and only if the greatest eigenvalue of $A$, denoted by $r(A)$, is strictly less than $1$. This follows from the behaviour of $A^mx^{(0)}$. In general we cannot obtain exact solution of $x=Ax+b$ by means of
Jacobi iterations.

In the case of $\rmax$, the situation is determined by the behaviour of $A^mx^{(0)}$
which differs from the case of the usual nonnegative algebra. However, this behaviour can be also analysed in terms of $r(A)$, the greatest eigenvalue
in terms of max-plus algebra (that is, with respect to the max-plus eigenproblem
$A\odot x=\lambda\odot x$).
 Namely, $A^m x^{(0)}\to\0$ and hence the iterations converge if $r(A)<\1$. Moreover $A^*=(I\oplus A\oplus\ldots\oplus A^{n-1})$ and hence the iterations yield {\bf exact} solution
to Bellman equation after a {\bf finite} number of steps.
To the contrary, $A^m x^{(0)}\to+\infty$ and hence the iterations diverge if $r(A)>\1$. See, for instance,~\cite{Car-71} for more details.
On the boundary $r(A)=1$, the powers $A^m$ reach a periodic regime
after a finite number of steps. Hence $A^*b\oplus A^m x^{(0)}$ also becomes
periodic, in general. If the period of $A^mx^{(0)}$ is one, that is, if this
sequence stabilizes, then the method converges to a general solution of
$x=Ax\oplus b$ described as a superposition of $A^*b$ and an eigenvector of $A$~\cite{BSS-11,Kri-06}. The vector $A^*b$ may dominate, in which case the
method converges to $A^*b$ as ``expected''. However, the period of $A^*b\oplus A^m x^{(0)}$ may be more than one, in which case the Jacobi iterations do not yield
any solution of $x=Ax\oplus b$. See~\cite{But:10} for more information on the
behaviour of max-plus matrix powers and the max-plus spectral theory.

In a more elaborate scheme of Gauss-Seidel iterations we can also
use the previously found coordinates of $X^{(k)}$. In this case matrix
$A$ is written as $L\oplus U$ where $L$ is the strictly lower triangular
part of $A$, and $U$ is the upper triangular part with the diagonal.
The iterations are written as
\begin{equation}
\label{e:seidel}
X^{(k)}=LX^{(k)}\oplus UX^{(k-1)}\oplus B=L^*UX^{(k-1)}\oplus L^*B
\end{equation}
Note that the transformation on the right hand side is unambiguous since
$L$ is strictly lower triangular and $L^*$ is uniquely defined as
$I\oplus L\oplus\ldots\oplus L^{n-1}$ (where $n$ is the dimension of $A$).
In other words, we just apply the forward substitution.
Iterating expressions
~\eqref{e:seidel} for all $k$ up to $m$ we obtain
\begin{equation}
\label{e:seid-unfold}
X^{(m)}=(L^*U)^m X^{(0)}\oplus \bigoplus_{i=0}^{m-1} (L^*U)^iL^*B
\end{equation}
The right hand side reminds of the formula
$(L\oplus U)^*=(L^*U)^*L^*$, see~\eqref{e:conway},
so it is natural to expect that these iterations converge to
$A^*B$ with a good choice of $X^{(0)}$.
The result crucially depends on the behaviour of $(L^*U)^m X^{(0)}$.
The algorithm can be written as follows (we assume again that $B$ is a column vector).

\begin{algorithm}
\label{a:seidel}
Gauss-Seidel iterations
\end{algorithm}

\sssn {\bf Input:} $n\times n$ matrix $A$ with entries $a(i,j)$;\\
$n\times 1$ column vectors $b$ and $x$

\msn situation$=$'proceed'\\
{\bf while}  situation$==$'proceed'\\
{\bf for} $i=1: n$\\
$y(i)= a(i,i: n)\odot x(i: n)\oplus b(i)$\\
{\bf end}\\
{\bf for} $i=2: n$\\
$x(i)=a(i,1: i-1)\odot x(1: i-1)$\\
{\bf end}\\
situation$=${\em newsituation}(...)\\
{\bf if} situation$==$'no convergence'\\
{\bf disp}('Gauss-Seidel iterations did not converge')\\
{\em exit}\\
{\bf end}\\
{\bf if} situation$==$'convergence'\\
{\bf disp}('Gauss-Seidel iterations converged')\\
{\em exit}\\
{\bf end}\\
{\bf end}

\msn {\bf Output:} situation, $x$.

\ssn It is plausible to expect that the behaviour of Gauss-Seidel scheme in the case of max-plus algebra and nonnegative linear algebra is analogous to the case of
Jacobi iterations.

\subsection{Software implementation of universal algorithms}
\label{ss:impl}

Software implementations for universal semiring algorithms cannot be
as efficient as hardware ones (with respect to the computation speed)
but they are much more flexible. Program modules can deal with abstract (and
variable) operations and data types. Concrete values for these
operations and data types can be defined by the corresponding
input data. In this case concrete operations and data types are generated
by means of additional program modules. For programs written in
this manner it is convenient to use special techniques of the
so-called object oriented (and functional) design, see,
for instance,~\cite{Lor:93,Pohl:97,SL:94}. Fortunately, powerful tools supporting the
object-oriented software design have recently appeared including compilers
for real and convenient programming languages (for instance, $C++$ and Java) and modern computer algebra systems. Recently, this type of programming technique has been dubbed
generic programming (see, for instance,~\cite{Pohl:97}).

{\em $C++$ implementation} Using templates and ob\-jec\-tive
ori\-en\-ted programming, Churkin and Ser\-geev~\cite{CS-07} created a Visual $C++$ application demonstrating how the
universal algorithms calculate matrix closures $A^*$ and solve Bellman equations
$x=Ax\oplus b$ in various semirings. The program can also compute the usual system $Ax=b$
in the usual arithmetic by transforming it to the ``Bellman'' form.
Before pressing ``Solve'', the user has to choose a semiring, a problem and an algorithm to use. Then the initial data are written into the
matrix (for the sake of visualization the dimension of a matrix is no more than $10\times 10$).
The result may appear as a matrix or as a vector depending on the problem to
solve. The object-oriented approach allows to implement various
semirings as objects with various definitions of basic operations, while keeping
the algorithm code unique and concise.

{\em Examples of the semirings.} The choice of semiring determines the
object used by the algorithm, that is, the concrete realization
of that algorithm. The following semirings have been realized:
\begin{itemize}
\item[1)] $\oplus=+$ and $\otimes=\times$: the usual arithmetic over reals;
\item[2)] $\oplus=\max$ and $\otimes=+$: max-plus arithmetic over $\R\cup\{-\infty\}$;
\item[3)] $\oplus=\min$ and $\otimes=+$: min-plus arithmetic over $\R\cup\{+\infty\}$;
\item[4)] $\oplus=\max$ and $\otimes=\times$: max-times arithmetic over
nonnegative numbers;
\item[5)] $\oplus=\max$ and $\otimes=\min$: max-min arithmetic over
a real interval $[a,b]$ (the ends $a$ and $b$ can be chosen by the user);
\item[6)] $\oplus=$OR and $\otimes=$AND: Boolean logic over the two-element
set $\{0,1\}$.
\end{itemize}

{\em Algorithms.} The user can select the following basic methods:
\begin{itemize}
\item[1)] {\bf Gaussian elimination scheme}, including the
universal realizations of escalator method (Algorithm~\ref{a:bordering} ),
Floyd-Warshall (Algorithm~\ref{a:eliminate},
Yershov's algorithm (based on a prototype from~\cite{FF} Ch. 2),
and the universal algorithm of Rote~\cite{Rot-85};
\item[2)] {\bf Methods for Toeplitz systems} including the universal
realizations of Durbin's and Levinson's schemes (Algorithms~\ref{a:YW} and
~\ref{a:levinson});
\item[3.] {\bf LDM decomposition} (Algorithm~\ref{a:ldm}) and its adaptations
to the symmetric case (Algorithm~\ref{a:ldl}), band matrices (Algorithm~\ref{a:ldmband}), Hessenberg and tridiagonal matrices.
\item[4)] {\bf Iteration schemes} of Jacobi and Gauss-Seidel. As mentioned above,
these schemes are not truly universal since the stopping criterion is different for
the usual arithmetics and idempotent semirings.
\end{itemize}

{\em Types of matrices.} The user may choose to work with general matrices, or
with a matrix of special structure, for instance, symmetric, symmetric Toeplitz, band,
Hessenberg or tridiagonal.

{\em Visualization.} In the case of idempotent semiring, the matrix can be visualized as a weighted digraph. After performing the calculations, the user
may wish to find an optimal path between a given pair of nodes, or to display an optimal paths tree. These problems can be solved using parental links like
in the case of the classical Floyd-Warshall meth\-od computing all optimal paths,
see, for instance, \cite{Sedge}. In our case, the mechanism of parental links can be implemented directly in the class describing an idempotent arithmetic.

{\em Other arithmetics and interval extensions.} It is also possible
to realize various types of arithmetics as data types and combine this with
the semiring selection. Moreover, all implemented semirings can be extended to their
interval versions. Such possibilities were not realized in the program of
Churkin and Sergeev~\cite{CS-07}, being postponed to the next version.
The list of such arithmetics includes integers, and  fractional arithmetics with the
use of chain fractions and controlled precision.

{\em MATLAB realization.} The whole work (except for visualization tools)
has been duplicated in MATLAB~\cite{CS-07}, which also allows for a kind of object-oriented programming. Obviously, the universal algorithms written in MATLAB are very close to those described in the present paper.

{\em Future prospects}. High-level tools, such as STL
\cite{Pohl:97,SL:94}, possess both obvious advantages and some disadvantages and must be used with caution.
It seems that it is natural to obtain an implementation of the correspondence
principle approach to scientific calculations in the form of a
powerful software system based on a collection of universal
algorithms. This approach should ensure a working time reduction for
programmers and users because of the software unification.
The arbitrary necessary accuracy and safety of numeric calculations can be ensured
as well.

The system has to contain several levels (including programmer and
user levels) and many modules.

Roughly speaking, it must be divided into three parts. The first part
contains modules that implement domain
modules (finite representations of
basic mathematical objects). The second part implements universal
(invariant) calculation methods. The third part contains modules
implementing model dependent algorithms. These modules may be
used in user programs written in $C^{++}$, Java, Maple, Matlab etc.

The system has to contain the following modules:

\medskip
\begin{itemize}
\item[---] Domain modules:
\begin{itemize}
\item infinite precision integers;
\item rational numbers;
\item finite precision rational numbers (see \cite{Ser-toep});
\item finite precision complex rational numbers;
\item fixed- and floating-slash rational numbers;
\item complex rational numbers;
\item arbitrary precision floating-point real numbers;
\item arbitrary precision complex numbers;
\item $p$-adic numbers;
\item interval numbers;
\item ring of polynomials over different rings;
\item idempotent semirings;
\item interval idempotent semirings;
\item and others.
\end{itemize}
\item[---] Algorithms:
\begin{itemize}
\item linear algebra;
\item numerical integration;
\item roots of polynomials;
\item spline interpolations and approximations;
\item rational and polynomial interpolations and approximations;
\item special functions calculation;
\item differential equations;
\item optimization and optimal control;
\item idempotent functional analysis;
\item and others.
\end{itemize}
\end{itemize}

This software system may be especially useful for designers
of algorithms, software engineers, students and mathematicians.

\if{
Note that there are some software systems oriented to calculations with idempotent semirings like $\rmax$; see, for instance~\cite{Scilab}. However these systems do not support universal algorithms.
}\fi

\section*{Acknowledgement}

The authors are grateful to the anonymous referees for a number of
important corrections in the paper.


\end{document}